# Paraconsistent second order arithmetic $\mathbf{Z}_2^{\#}$ based on the paraconsistent logic $\overline{\mathbf{LP}}_{\omega}^{\#}$ with infinite hierarchy levels of contradiction. Schinzel's conjecture H within paraconsistent field $\mathbb{C}_{\mathbf{inc}}^{\#}$.

Jaykov Foukzon
Israel Institute of Technology
jaykovfoukzon@list.ru

**Abstract**: In this paper paraconsistent second order arithmetic $\mathbf{Z}_2^{\#}$ with unrestrictedcomprehension scheme is proposed.We outline the development of certain portions of paraconsistent mathematics within paraconsistent second order arithmetic $\mathbf{Z}_2^{\#}$. In particular we defined infinite hierarchy Berry's and Richard's inconsistent numbers as elements of the paraconsistent field $\mathbb{R}^{\#}$.

**Contents**

References.

# I.Introduction.

Let be $\mathbf{Z}_2^*$ second order arithmetic [3]-[5] with second order language $L_2$ and with unrestricted comprehension scheme:

$$\exists X \forall n (n \in X \leftrightarrow \varphi(n,X)) \tag{1.1}$$

It is known that second order arithmetic $\mathbf{Z}_2^*$ is inconsistent from the well known standard construction named as Berry's and Richard's inconsistent numbers. Suppose that $\mathcal{F}(n,X) \in L_2$ is a well-formed formula of second-order arithmetic $\mathbf{Z}_2^*$, i.e. formula which is arithmetical, which has one free set variable $X$ and one free individual variable $n$. Suppose that $g(\exists X \mathcal{F}(x,X)) \leq \mathbf{k}$, where $g(\exists X \mathcal{F}(x,X))$ is a corresponding Gödel number. Let be $A_{\mathbf{k}}, \mathbf{k} \in \mathbb{N}$ the set of all positive consistent integers $\bar{n}$ which can be defined under corresponding well-formed formula $\mathcal{F}_{\bar{n}}(x,X)$, i.e. $\exists X_{\bar{n}} \forall m [\mathcal{F}_{\bar{n}}(m,X_{\bar{n}}) \to m = \bar{n}]$, hence $\bar{n} \in A_{\mathbf{k}} \leftrightarrow \exists X \mathcal{F}_{\bar{n}}(\bar{n},X_{\bar{n}})$.

Thus

$$\begin{gathered} \forall n [n \in A_{\mathbf{k}} \leftrightarrow \exists X_n \mathcal{F}_n(n,X_n)], \\ g(\exists X_n \mathcal{F}_n(x,X_n)) \leq \mathbf{k}, \end{gathered} \tag{1.2}$$

where $g(\exists X \mathcal{F}(x,X))$ is a corresponding Gödel number.Since there are only finitely many of these $\bar{n}$, there must be a smallest positive integer $n_{\mathbf{k}} \in \mathbb{N} \backslash A_{\mathbf{k}}$ that does not belong to $A_{\mathbf{k}}$. But we just defined $n_{\mathbf{k}}$ in under corresponding well-formed formula

$$\begin{gathered} n_{\mathbf{k}} \in A_{\mathbf{k}} \leftrightarrow \breve{\mathcal{F}}_{n_{\mathbf{k}}}(n_{\mathbf{k}},A_{\mathbf{k}}), \\ \breve{\mathcal{F}}_{n_{\mathbf{k}}}(n_{\mathbf{k}},A_{\mathbf{k}}) \leftrightarrow n_{\mathbf{k}} = \min_{n \in \mathbb{N}} (\mathbb{N} \backslash A_{\mathbf{k}}). \end{gathered} \tag{1.3}$$

Hence for a sufficiently Large $\mathbf{k}$ such that: $g\left(\breve{\mathcal{F}}(n_{\mathbf{k}},A_{\mathbf{k}})\right) \leq \mathbf{k}$ one obtain the contradiction: $(n_{\mathbf{k}} \in A_{\mathbf{k}}) \wedge (n_{\mathbf{k}} \notin A_{\mathbf{k}})$.

Within $\mathbf{Z}_2^*$, a consistent real number $x = \langle q_n : n \in \mathbb{N} \rangle \in \mathbb{R}$ is defined to be a Cauchy consistent sequence $\langle q_n | n \in \mathbb{N}\rangle$ of rational numbers, i.e., a consistent sequence of rational numbers $x = \langle q_n : n \in \mathbb{N}, q_n \in \mathbb{Q}\rangle$ such that

$$\forall \varepsilon(\varepsilon \in \mathbb{Q})(\varepsilon > 0 \to \exists m \forall n(m < n \to |q_m - q_n| < \varepsilon)), \qquad (1.4)$$

see **Definition 2.2.9**.

Let be $q_n \in \mathbb{Q}$ rational number with corresponding decimal representation $q_n = \{0, q_n(1)q_n(2)\ldots q_n(n)\}, q_n(i) = 0,1,2,\ldots,9, i \leq n$, $x_k = \langle q_n^k\rangle = \langle q_n^k : n \in \mathbb{N}, q_n^k = \{0, q_n^k(1)q_n^k(2)\ldots q_n^k(n)\}\rangle \in \mathbb{R}$ is a consistent real number which can be defined under corresponding well-formed formula (of second-order arithmetic $\mathbf{Z}_2$) $\mathcal{F}_k(x)$, i.e. $\forall q(q \in \mathbb{Q})[q \in \langle q_n^k\rangle \leftrightarrow \exists X \mathcal{F}_k(q, X)]$. We denote real number $x_k$ as $k$-th Richard's real number.

Let us consider Richard's real number $\Re_p = \langle \Re_n^p : n \in \mathbb{N}\rangle$ such that

$$\Re_n^p = 1 \leftrightarrow q_n^n(n) \neq 1,$$
$$\Re_n^p = 0 \leftrightarrow q_n^n(n) = 1. \qquad (1.5)$$

Suppose that $q_p^p(p) \neq 1$, hence $\Re_p^p(p) = 1$. Thus $\Re_p^p(p) \neq q_p^p(p) \to \Re_p \neq x_p$.
Suppose that $q_p^p(p) = 1$, hence $\Re_p^p(p) = 0$. Thus $\Re_p^p(p) \neq q_p^p(p) \to \Re_p \neq x_p$.

Hence for any Richard's real number $x_k$ one obtain the contradiction $x_k \neq \Re_p^p(p)$. Thus the classical logical antinomy known as Richard-Berry paradox is combined with plausible assumptions formalizing certain sentences, to show that formalization of language leads to contradictions which trivialize the system $\mathbf{Z}_2^*$. In this paper paraconsistent second order arithmetic $\mathbf{Z}_2^\#$ with unrestricted comprehension scheme is proposed. We outline the development of certain portions of paraconsistent mathematics within paraconsistent second order arithmetic $\mathbf{Z}_2^\#$. In particular we defined infinite hierarchy Berry's and Richard's inconsistent numbers as elements of the paraconsistent field $\mathbb{R}^\#$.

# II.Consistent second order arithmetic $\mathbf{Z}_2$.

## II.1.System $\mathbf{Z}_2$.

In this section we briefly define $\mathbf{Z}_2$, the well known formal system of second order consistent arithmetic. For more detailed information concerning this system see [3]-[5].

The language $L_2$ of second order consistent arithmetic is a two-sorted language. This means that there are two distinct sorts of variables which are intended to range over two different kinds of object.

(**1**) **Variables**:

(**1.1**) Variables of the *first sort*: are known as consistent number variables, are denoted by $i, j, k, m, n, \ldots,$ and are intended to range over the set $\omega = \{0, 1, 2, \ldots\}$ of all consistent natural numbers.

(**1.2**) Variables of the second sort are known as consistent set variables, are denoted by $X, Y, Z, \ldots,$ and are intended to range over all subsets of $\omega$.

The terms and formulas of the language of second order consistent arithmetic are as follows:

(**2**) Numerical terms are number variables, the constant symbols $0$ and $1,$ and $\mathbf{t}_1 + \mathbf{t}_2$ and $\mathbf{t}_1 \times \mathbf{t}_2$ whenever $t_1$ and $t_2$ are numerical terms.

Here $(\bullet + \bullet)$ and $(\bullet \times \bullet)$ are binary operation symbols intended to denote addition and multiplication of natural numbers. (Numerical terms are intended to denote natural numbers.)

(**3**) Atomic formulas are:

(**3.1**) $\mathbf{t}_1 = \mathbf{t}_2, \mathbf{t}_1 < \mathbf{t}_2,$

(**3.2**) $\mathbf{t}_1 \in X$ where $\mathbf{t}_1$ and $\mathbf{t}_2$ are numerical terms and $X$ is any set variable.

(The intended meanings of these respective atomic formulas are that $\mathbf{t}_1$ equals $\mathbf{t}_2, \mathbf{t}_1$ is less than $\mathbf{t}_2,$ and $\mathbf{t}_1$ is an element of $X.$)

(**4**) Formulas are built up from:

(**4.1**) atomic formulas by means of propositional connectives $\wedge, \vee, \neg, \rightarrow, \leftrightarrow$ (and, or, not, implies, if and only if),

(**4.2**) *consistent number* quantifiers $\forall n, \exists n$ (for all $n,$ there exists $n$),

(**4.3**) *consistent set* quantifiers $\forall X, \exists X$ (for all $X,$ there exists $X$).

(**5**) A sentence is a formula with no free variables.

**Definition 2.1.1**. (language of second order consistent arithmetic). $L_2$ is defined to be the language of second order consistent arithmetic as described above.

In writing terms and formulas of $L_2,$ we shall use parentheses and brackets to indicate grouping, as is customary in mathematical logic textbooks. We shall also use some obvious abbreviations. For instance, $2 + 2 = 4$ stands for $(1 + 1) + (1 + 1) = ((1 + 1) + 1) + 1,$ $(m + n)^2 \notin X$ stands for $\neg((m + n) \cdot (m + n) \in X), \mathbf{s} \leq \mathbf{t}$ stands for $\mathbf{s} < \mathbf{t} \vee \mathbf{s} = \mathbf{t},$ and $\varphi \wedge \psi \wedge \theta$ stands for $(\varphi \wedge \psi) \wedge \theta.$

The semantics of the language $L_2$ are given by the following definition.

**Definition 2.1.2**. ($L_2$-structures). A model for $L_2,$ also called a structure for $L_2$ or an $L_2$-structure, is an ordered $7$-tuple:

$$\mathbf{M} = (|\mathbf{M}|, \mathbf{S_M}, (\cdot +_\mathbf{M} \cdot), (\cdot \times_\mathbf{M} \cdot), 0_\mathbf{M}, 1_\mathbf{M}, (\cdot <_\mathbf{M} \cdot)), \qquad (2.1.1)$$

where $|\mathbf{M}|$ is a set which serves as the range of the number variables, $\mathbf{S_M}$ is a set of subsets of $|\mathbf{M}|$ serving as the range of the set variables, $+_\mathbf{M}$ and $\times_\mathbf{M}$ are binary operations on $|\mathbf{M}|$, $0_\mathbf{M}$ and $1_\mathbf{M}$ are distinguished elements of $|\mathbf{M}|$, and $<_\mathbf{M}$ is a binary relation on $|\mathbf{M}|$. We always assume that the sets $|\mathbf{M}|$ and $\mathbf{S_M}$ are disjoint and nonempty. Formulas of $L_2$ are interpreted in $\mathbf{M}$ in the obvious way.

In discussing a particular model $\mathbf{M}$ as above, it is useful to consider formulas with parameters from $|\mathbf{M}| \cup \mathbf{S_M}$. We make the following slightly more general definition.

**Definition 2.1.3**. (parameters). Let $\mathbf{B}$ be any subset of $|\mathbf{M}| \cup \mathbf{S_M}$. By a formula with parameters from $\mathbf{B}$ we mean a formula of the extended language $L_2(\mathbf{B})$. Here $L_2(\mathbf{B})$ consists of $L_2$augmented by new constant symbols corresponding to the elements of $\mathbf{B}$. By a sentence with parameters from $\mathbf{B}$ we mean a sentence of $L_2(\mathbf{B})$, i.e.,a formula of $L_2(\mathbf{B})$ which has no free variables.

In the language $L_2(|\mathbf{M}| \cup \mathbf{S_M})$, constant symbols corresponding to elements of $\mathbf{S_M}$ (respectively $|\mathbf{M}|$) are treated syntactically as unquantified set variables (respectively unquantified number variables). Sentences and formulas with parameters from $|\mathbf{M}| \cup \mathbf{S_M}$ are interpreted in $\mathbf{M}$ in the obvious way.

**Definition 2.1.4**. A set $\mathbf{A} \subseteq |\mathbf{M}|$ is said to be definable over $\mathbf{M}$ allowing parameters from $\mathbf{B}$ if there exists a formula $\varphi(n)$ with parameters from $\mathbf{B}$ and no free variables other than $n$ such that $A = \{a \in |\mathbf{M}| : \mathbf{M} \models \varphi(a)\}$.Here $\mathbf{M} \models \varphi(a)$ means that $\mathbf{M}$ satisfies $\varphi(a)$,i.e., $\varphi(a)$ is true in $\mathbf{M}$.

**Definition 2.1.5**.The *intended model* for $L_2$ is of course the model $(\omega, \mathbf{P}(\omega), +, \times, 0, 1, <)$ where $\omega$ is the set of natural numbers, $\mathbf{P}(\omega)$ is the set of all subsets of $\omega$, and $+, \cdot, 0, 1, <$ are as usual.

By an $\omega$-model we mean an $L_2$-structure of the form $(\omega, \mathbf{S}, +, \cdot, 0, 1, <)$ where $\emptyset \neq \mathbf{S} \subseteq \mathbf{P}(\omega)$. Thus an $\omega$-model differs from the intended model only by having a possibly smaller collection $\mathbf{S}$ of sets to serve as the range of the set variables. We sometimes speak of the $\omega$-model $\mathbf{S}$ when we really mean the $\omega$-model $(\omega, \mathbf{S}, +, \cdot, 0, 1, <)$.

**Definition 2.1.6**. (second order arithmetic $\mathbf{Z}_2$). The axioms of second order arithmetic $\mathbf{Z}_2$ consist of the universal closures of the following $L_2$-formulas:

(**i**) basic axioms:

(**i.1**) $n + 1 \neq 0$,

(**i.2**) $m + 1 = n + 1 \to m = n$,

(**i.3**) $m + 0 = m$,

(**i.4**) $m + (n + 1) = (m + n) + 1$,

(**i.5**) $m \times 0 = 0$,

(**i.6**) $m \times (n + 1) = (m \times n) + m$,

(**i.7**) $\neg(m < 0)$,
(**i.8**) $m < n+1 \leftrightarrow (m < n \vee m = n)$.
(**ii**) induction axiom:
$(0 \in X \wedge \forall n(n \in X \to n+1 \in X)) \ \to \forall n(n \in X)$
(**iii**) comprehension scheme:
$\exists X \forall n(n \in X \leftrightarrow \varphi(n))$
where $\varphi(n)$ is any formula of $L_2$ in which $X$ does not occur freely.

## II.2.Consistent Mathematics Within $\mathbf{Z}_2$.

We now outline the development of certain portions of ordinary mathematics within $\mathbf{Z}_2$.

**Definition 2.2.1.** If $X$ and $Y$ are set variables, we use $X = Y$ and $X \subseteq Y$ as abbreviations for the formulas $\forall n(n \in X \leftrightarrow n \in Y)$ and $\forall n(n \in X \to n \in Y)$ respectively.

**Definition 2.2.2.** Within $\mathbf{Z}_2$, we define $\mathbb{N}$ to be the unique set $X$ such that $\forall n(n \in X)$.

**Definition 2.2.3.** For $X, Y \subseteq \mathbb{N}$, a consistent function $f : X \to Y$ is defined to be a consistent set $f \subseteq X \times Y$ such that for all $m \in X$ there is exactly one $n \in Y$ such that $(m,n) \in f$. For $m \in X$, $f(m)$ is defined to be the unique $n$ such that $(m,n) \in f$. The usual properties of such functions can be proved in $\mathbf{Z}_2$.

**Definition 2.2.4.** (consistent primitive recursion). This means that, given $f : X \to Y$ and $g : \mathbb{N} \times X \times Y \to Y$, there is a unique $h : \mathbb{N} \times X \to Y$ defined by
$h(0,m) = f(m)$,
$h(n+1,m) = g(n,m,h(n,m))$ for all $n \in \mathbb{N}$ and $m \in X$.
The existence of $h$ is proved by arithmetical comprehension, and the uniqueness of $h$ is proved by arithmetical induction.
In particular, we have the exponential function $\exp(m,n) = m^n$, defined by $m^0 = 1, m^{n+1} = m^n \times m$ for all $m,n \in \mathbb{N}$. The usual properties of the exponential function can be proved in $\mathbf{Z}_2$.
The consistent natural number system is essentially already given to us by the language $L_2$ and axioms of $\mathbf{Z}_2$. Thus, within $\mathbf{Z}_2$, a consistent natural number is defined to be an element of $\mathbb{N}$, and the natural number system is defined to be the structure $\mathbb{N}, +_{\mathbb{N}}, \times_{\mathbb{N}}, 0_{\mathbb{N}}, 1_{\mathbb{N}}, <_{\mathbb{N}}, =_{\mathbb{N}}$, where $+_{\mathbb{N}} : \mathbb{N} \times \mathbb{N} \to \mathbb{N}$ is defined by $m +_{\mathbb{N}} n = m+n$, etc. Thus for instance $+_{\mathbb{N}}$ is the set of triples $((m,n),k) \in (\mathbb{N} \times \mathbb{N}) \times \mathbb{N}$ such that $m+n = k$. The existence of this set follows from the arithmetical comprehension.

In a standard manner, we can define within $\mathbf{Z}_2$ the set $\mathbb{Z}$ of consistent integers and the set of consistent rational numbers: $\mathbb{Q}$.

**Definition 2.2.5**.(consistent rational numbers $\mathbb{Q}$) Let $\mathbb{Z}^+ = \{a \in \mathbb{Z} : 0 <_{\mathbb{Z}} a\}$ be the set of positive consistent integers,and let $\equiv_{\mathbb{Q}}$ be the equivalence relation on $\mathbb{Z} \times \mathbb{Z}^+$ defined by $(a,b) \equiv_{\mathbb{Q}} (c,d)$ if and only if $a \times_{\mathbb{Z}} d = b \times_{\mathbb{Z}} c$. Then $\mathbb{Q}$ is defined to be the set of all $(a,b) \in \mathbb{Z} \times \mathbb{Z}^+$ such that $(a,b)$ is the $<_{\mathbb{N}}$ -minimum element of its $\equiv_{\mathbb{Q}}$ -equivalence class. Operations $+_{\mathbb{Q}}, -_{\mathbb{Q}}, \times_{\mathbb{Q}}$ on $\mathbb{Q}$ are defined by:
$(a,b) +_{\mathbb{Q}} (c,d) \equiv_{\mathbb{Q}} (a \times_{\mathbb{Z}} d +_{\mathbb{Z}} b \times_{\mathbb{Z}} c, b \times_{\mathbb{Z}} d), -_{\mathbb{Q}} (a,b) \equiv_{\mathbb{Q}} (-_{\mathbb{Z}} a, b)$, and
$(a,b) \times_{\mathbb{Q}} (c,d) \equiv Q(a \times_{\mathbb{Z}} c, b \times_{\mathbb{Z}} d)$. We let $0_{\mathbb{Q}} \equiv_{\mathbb{Q}} (0_{\mathbb{Z}}, 1_{\mathbb{Z}})$ and $1_{\mathbb{Q}} \equiv_{\mathbb{Q}} (1_{\mathbb{Z}}, 1_{\mathbb{Z}})$, and we define a binary relation $<_{\mathbb{Q}}$ on $\mathbb{Q}$ by letting $(a,b) <_{\mathbb{Q}} (c,d)$ if and only if $a \times_{\mathbb{Z}} d <_{\mathbb{Z}} b \times_{\mathbb{Z}} c$. Finally $=_{\mathbb{Q}}$ is the identity relation on $\mathbb{Q}$. We can then prove within $\mathbf{Z}_2$ that the rational number system $\mathbb{Q}, +_{\mathbb{Q}}, -_{\mathbb{Q}}, \times_{\mathbb{Q}}, 0_{\mathbb{Q}}, 1_{\mathbb{Q}}, <_{\mathbb{Q}}, =_{\mathbb{Q}}$ has the usual properties ofan ordered field, etc.

We make the usual identifications whereby $\mathbb{N}$ is regarded as a subset of $\mathbb{Z}$ and $\mathbb{Z}$ is regarded as a subset of $\mathbb{Q}$. Namely m $\in \mathbb{N}$ is identified with $(m,0) \in \mathbb{Z}$, and $a \in \mathbb{Z}$ is identified with $(a, 1_{\mathbb{Z}}) \in \mathbb{Q}$. We use $+$ ambiguously to denote $+_{\mathbb{N}}, +_{\mathbb{Z}}$, or $+_{\mathbb{Q}}$ and similarly for $-, \times, 0, 1, <$. For $q, r \in \mathbb{Q}$ we write $q - r = q + (-r)$, and if $r \neq 0, q/r =$ the unique $q' \in \mathbb{Q}$ such that $q = q' \times r$. The function $\exp(q,a) = q^a$ for $q \in \mathbb{Q} \backslash \{0\}$ and $a \in \mathbb{Z}$ is obtained by primitive recursion in the obvious way.

**Definition 2.2.6**.The absolute value function $|\cdot| : \mathbb{Q} \to \mathbb{Q}$ is defined by $|q| = q$ if $q \geq 0, -q$ otherwise.

**Definition 2.2.7**. An consistent sequence of rational numbers is defined to be a consistent function $f : \mathbb{N} \to \mathbb{Q}$.
We denote such a sequence as $\langle q_n : n \in \mathbb{N} \rangle$, or simply $\langle q \rangle_n$, where $q_n = f(n)$.

**Definition 2.2.8**. A double consistent sequence of rational numbers is a consistent function $f : \mathbb{N} \times \mathbb{N} \to \mathbb{Q}$, denoted $\langle q_{mn} : m,n \in \mathbb{N} \rangle$ or simply $\langle q_{mn} \rangle$, where $q_{mn} = f(m,n)$.

**Definition 2.2.9**.(consistent real numbers). Within $\mathbf{Z}_2$, a consistent real number is defined to be a Cauchy consistent sequence of rational numbers, i.e., a consistent sequence of rational numbers $x = \langle q_n : n \in \mathbb{N} \rangle$ such that

$$\forall \varepsilon (\varepsilon \in \mathbb{Q})(\varepsilon > 0 \to \exists m \forall n (m < n \to |q_m - q_n| < \varepsilon)). \qquad (2.2.1)$$

**Definition 2.2.10**. If $x = q_n$ and $y = q'_n$ are consistent real numbers, we write $x =_{\mathbb{R}} y$ to mean that $\lim_n |q_n - q'_n| = 0$, i.e.,

$$\forall \varepsilon(\varepsilon \in \mathbb{Q})(\varepsilon > 0 \to \exists m \forall n(m < n \to |q_n - q_n'| < \varepsilon)), \qquad (2.2.2)$$

and we write $x <_{\mathbb{R}} y$ to mean that

$$\exists \varepsilon(\varepsilon > 0 \wedge \exists m \forall n(m < n \to q_n + \varepsilon < q_n')). \qquad (2.2.3)$$

Also $x +_{\mathbb{R}} y = \langle q_n + q_n' \rangle, x \times_{\mathbb{R}} y = \langle q_n \times q_n' \rangle, -_{\mathbb{R}} x = \langle -q_n \rangle, 0_{\mathbb{R}} = \langle 0 \rangle, 1_{\mathbb{R}} = \langle 1 \rangle$.
We use $\mathbb{R}$ to denote the set of all *consistent real numbers*. Thus $x \in \mathbb{R}$ means that $x$ is a *consistent real number.* (Formally, we cannot speak of the set $\mathbb{R}$ within the language of second order arithmetic, since it is a set of sets.)
We shall usually omit the subscript $\mathbb{R}$ in $+_{\mathbb{R}}, -_{\mathbb{R}}, \times_{\mathbb{R}}, 0_{\mathbb{R}}, 1_{\mathbb{R}}, <_{\mathbb{R}}, =_{\mathbb{R}}$.
Thus the *consistent real number system* consists of $\mathbb{R}, +, -, \times, 0, 1, <, =$. We shall sometimes identify a consistent rational number $q \in \mathbb{Q}$ with the corresponding consistent real number $x_q = \langle q \rangle$.
Within $\mathbf{Z}_2$ one can prove that the real number system has the usual properties of an *consistent Archimedean ordered field*, etc. The complex consistent numbers can be introduced as usual as pairs of real numbers.
Within $\mathbf{Z}_2$,it is straightforward to carry out the proofs of all the basic results in real and complex linear and polynomial algebra. For example, the fundamental theorem of algebra can be proved in $\mathbf{Z}_2$.
**Definition 2.2.11**. A consistent sequence of real numbers is defined to be a double consistent sequence of rational numbers $\langle q_{mn} : m, n \in \mathbb{N} \rangle$ such that for each $m, \langle q_{mn} : n \in \mathbb{N} \rangle$ is a consistent real number. Such a sequence of real numbers is denoted $\langle x_m : m \in \mathbb{N} \rangle$, where $x_m = \langle q_{mn} : n \in \mathbb{N} \rangle$. Within $\mathbf{Z}_2$ we can prove that every bounded consistent sequence of real numbers has a *consistent least upper bound.* This is a very useful completeness property of the consistent real number system.For instance, it implies that an infinite series of positive terms is convergent if and only if the finite partial sums are bounded.
We now turn of certain portions of consistent abstract algebra within $\mathbf{Z}_2$.
Because of the restriction to the language $L_2$ of second order arithmetic, we cannot expect to obtain a good general theory of arbitrary (countable and uncountable) algebraic structures. However, we can develop countable algebra,i.e., the theory of countable algebraic structures, within $\mathbf{Z}_2$.
**Definition 2.2.12**. A countable consistent commutative ring is defined within $\mathbf{Z}_2$ to be a consistent structure $\mathbf{R}, +_{\mathbf{R}}, -_{\mathbf{R}}, \times_{\mathbf{R}}, 0_{\mathbf{R}}, 1_{\mathbf{R}}$, where $\mathbf{R} \subseteq \mathbb{N}$, $+_{\mathbf{R}} : \mathbf{R} \times \mathbf{R} \to \mathbf{R}$, etc., and the usual commutative ring axioms are assumed.
(We include $0 \neq 1$among those axioms.) The subscript $\mathbf{R}$ is usually omitted.

An ideal in $\mathbf{R}$ is a set $I \subseteq \mathbf{R}$ such that $a \in I$ and $b \in I$ imply $a + b \in I; a \in I$ and $r \in \mathbf{R}$ imply $a \times r \in I,$ and $0 \in I$ and $1 \notin I.$ We define an equivalence relation $=_I$ on $\mathbf{R}$ by $r =_I s$ if and only if $r - s \in I.$ We let $\mathbf{R}/I$ be the set of $r \in \mathbf{R}$ such that $r$ is the $<_{\mathbb{N}}$ -minimum element of its equivalence class under $=_I$ .Thus $\mathbf{R}/I$ consists of one element of each $=_I$ -equivalence class of elements of $\mathbf{R}$. With the appropriate operations,$\mathbf{R}/I$ becomes a countable commutative ring, the quotient ring of $\mathbf{R}$ by $I$.
The ideal $I$ is said to be prime if $\mathbf{R}/I$ is an integral domain, and maximal if $\mathbf{R}/I$ is a field.
Next we indicate how some basic concepts and results of analysis and topology can be developed within $\mathbf{Z}_2$.
**Definition 2.2.13**.Within $\mathbf{Z}_2$,a complete separable consistent metric space is a nonempty set $A \subseteq \mathbb{N}$ together with a function $d : A \times A \to \mathbb{R}$ satisfying $d(a,a) = 0,\ d(a,b) = d(b,a) \geq 0,$ and $d(a,c) \leq d(a,b) + d(b,c)$ for all $a,b,c \in A.$ (Formally, $d$ is a consistent sequence of real numbers, indexed by $A \times A.$) We define a point of the complete separable metric space $\widehat{A}$ to be a sequence $x = \langle a_n : n \in \mathbb{N}\rangle,\ a_n \in A,$ satisfying

$$\forall\varepsilon(\varepsilon \in \mathbb{R})(\varepsilon > 0 \to \exists m \forall n(m < n \to d(a_m, a_n) < \varepsilon)). \tag{2.2.4}$$

The pseudometric $d$ is extended from $A$ to $\widehat{A}$ by

$$d(x,y) = \lim_{n\to\infty} d(a_n, b_n) \tag{2.2.5}$$

where $x = \langle a_n : n \in \mathbb{N}\rangle$ and $y = \langle b_n : n \in \mathbb{N}\rangle.$ We write $x = y$ if and only if $d(x,y) = 0.$ For example, $\mathbb{R} = \widehat{\mathbb{Q}}$ under the metric $d(q,q') = |q - q'|$.
**Definition 2.2.14**.(consistent continuous functions). Within $\mathbf{Z}_2$, if $\widehat{A}$ and $\widehat{B}$ are complete separable metric spaces,a consistent continuous function $\phi : \widehat{A} \to \widehat{B}$ is a set $\Phi \subseteq A \times \mathbb{Q}^+ \times B \times \mathbb{Q}^+$satisfying the following coherence conditions:

$$1.\,[(a,r,b,s)\in\Phi]\wedge[(a,r,b',s')\in\Phi]\dashrightarrow d(b,b')<s+s';$$

$$2.\,[(a,r,b,s)\in\Phi]\wedge[d(b,b')+s<s']\dashrightarrow(a,r,b',s')\in\Phi \quad (2.2.6)$$

$$3.\,[(a,r,b,s)\in\Phi]\wedge[d(a,a')+r'<r]\dashrightarrow(a',r',b,s)\in\Phi$$

# III.Paraconsistent second order arithmetic $\mathbf{Z}_2^{\#}$.

## III.1.Paraconsistent system $\mathbf{Z}_2^{\#}$.

In this section we define $\mathbf{Z}_2$, the formal system of second order paraconsistent arithmetic based on the paraconsistent logic $\overline{\mathbf{LP}}_{\omega}^{\#}$ [1] with infinite hierarchy levels of contradiction. For detailed information concerning paraconsistent logic $\mathbf{LP}_{\omega}^{\#}$ see [2].

**Definition 3.1.1**. For any binary inconsistent relation $(\cdot \sim_w \cdot)$ we define:

$$a{\sim_{w,(0)}}b \triangleq \left(a \sim_w b\right)^{(0)},\ldots,a{\sim_{w,(n)}}b \triangleq \left(a \sim_w b\right)^{(n)};$$

$$a{\sim_{w,[0]}}b \triangleq \left(a \sim_w b\right)^{[0]},\ldots,a{\sim_{w,[n]}}b \triangleq \left(a \sim_w b\right)^{[n]};n\in\mathbb{N}, \quad (3.1.1)$$

$$\mathbb{N}=\{0,1,2,\ldots\}.$$

Here (i) $\left(a \sim_w b\right)^{(n)}$ is defined as follows:

$$\left(a\sim_w b\right)^{(0)} \triangleq a\sim_{w,(0)} b \triangleq a\sim_w b,$$

$$\left(a\sim_w b\right)^{(n)} \triangleq \bigwedge_{i=1}^{n}\left(a\sim_w b\right)^{i}, \left(a\sim_w b\right)^{i+1} \triangleq \left(\left(a\sim_w b\right)^{i}\right)^{1}, \qquad (3.1.1.a)$$

$$\left(a\sim_w b\right)^{1} \triangleq \neg\left[\left(a\sim_w b\right)\wedge\neg\left(a\sim_w b\right)\right],$$

(ii) $\left(a\sim_w b\right)^{[n]}$ is defined as follows:

$$\left(a\sim_w b\right)^{[0]} \triangleq a\sim_{w,[0]} b \triangleq a\sim_w b,$$

$$\left(a\sim_w b\right)^{[n]} \triangleq \bigwedge_{i=1}^{n}\left(a\sim_w b\right)^{\{i\}}, \left(a\sim_w b\right)^{\{i+1\}} \triangleq \left(\left(a\sim_w b\right)^{\{i\}}\right)^{\{1\}}, \qquad (3.1.1.b)$$

$$\left(a\sim_w b\right)^{\{1\}} \triangleq \left[\left(a\sim_w b\right)\wedge\neg\left(a\sim_w b\right)\right]$$

**Definition 3.1.2**. In particular we define:

$$a =_{w,(0)} b \triangleq (a =_w b)^{(0)},\ldots,a =_{w,(n)} b \triangleq (a =_w b)^{(n)};$$

$$a =_{w,[0]} b \triangleq (a =_w b)^{[0]},\ldots,a =_{w,[n]} b \triangleq (a =_w b)^{[n]}; n \in \mathbb{N}.$$

$$(3.1.2)$$

$$a <_{w,(0)} b \triangleq (a <_w b)^{(0)},\ldots,a <_{w,(n)} b \triangleq (a <_w b)^{(n)};$$

$$a <_{w,[0]} b \triangleq (a <_w b)^{[0]},\ldots,a <_{w,[n]} b \triangleq (a <_w b)^{[n]}; n \in \mathbb{N}$$

The language $L_2^{\#}$ of second order paraconsistent arithmetic $\mathbf{Z}_2^{\#}$ is a two-sorted paraconsistent language. This means that there are two distinct sorts of variables which are intended to range over two different kinds of object.

(**1**) **Variables**:

(**1.1**) Variables of the *first sort*: are denoting *consistent and inconsistent* number variables, are denoted as in classical case by $i,j,k,m,n,\ldots,$ and are intended to range over the set $\mathbb{N}^{\#} \supsetneqq \mathbb{N} = \{0,1,2,\ldots\}$ of all *consistent and inconsistent natural numbers.*

(**1.2**) Variables of the *second sort*: are denoting *consistent and inconsistent* set variables,are denoted by $X,Y,Z,\ldots,$ and are intended to range over all subsets of $\mathbb{N}^{\#}$. The terms and formulas of the language $L_2^{\#}$ of second order paraconsistent arithmetic $\mathbf{Z}_2^{\#}$ are as follows:

(**2**) Numerical terms are number variables, the constant symbols $0_{\mathbf{s}}, 0_w, 0_{w,(i)}, 0_{w,[i]}, 1_{\mathbf{s}}, 1_w, 1_{w,(i)}, 1_{w,[i]}, i \in \mathbb{N}$ and $\mathbf{t}_1 + \mathbf{t}_2$ $\mathbf{t}_1 \times \mathbf{t}_2$ whenever $\mathbf{t}_1$ and $\mathbf{t}_2$ are numerical terms in general.

Here $(\cdot + \cdot)$ and $(\cdot \times \cdot)$ are binary operation symbols intended to denote addition and multiplication of *consistent and inconsistent* natural numbers. (Numerical terms are intended to denote *consistent and inconsistent* natural numbers.)

(**3**) Atomic formulas are:

(**3.1**) $\mathbf{t}_1 =_{\mathbf{s}} \mathbf{t}_2, \mathbf{t}_1 <_{\mathbf{s}} \mathbf{t}_2, \mathbf{t}_1 =_w \mathbf{t}_2, \mathbf{t}_1 <_w \mathbf{t}_2,$

(**3.2**) $\mathbf{t}_1 =_{w,(0)} \mathbf{t}_2, \mathbf{t}_1 =_{w,(1)} \mathbf{t}_2, \ldots, \mathbf{t}_1 =_{w,(n)} \mathbf{t}_2, n \in \mathbb{N},$

(**3.3**) $\mathbf{t}_1 =_{w,[0]} \mathbf{t}_2, \mathbf{t}_1 =_{w,[1]} \mathbf{t}_2, \ldots, \mathbf{t}_1 =_{w,[n]} \mathbf{t}_2, n \in \mathbb{N},$

(**3.4**) $\mathbf{t}_1 <_{w,(0)} \mathbf{t}_2, \mathbf{t}_1 <_{w,(1)} \mathbf{t}_2, \ldots, \mathbf{t}_1 <_{w,(n)} \mathbf{t}_2, n \in \mathbb{N},$

(**3.5**) $\mathbf{t}_1 <_{w,[0]} \mathbf{t}_2, \mathbf{t}_1 <_{w,[1]} \mathbf{t}_2, \ldots, \mathbf{t}_1 <_{w,[n]} \mathbf{t}_2, n \in \mathbb{N},$

(**3.6**) $\mathbf{t}_1 \in_{\mathbf{s}} X, \mathbf{t}_1 \in_w X,$

(**3.7**) $\mathbf{t}_1 \in_{w,(0)} \mathbf{t}_2, \mathbf{t}_1 \in_{w,(1)} \mathbf{t}_2, \ldots, \mathbf{t}_1 \in_{w,(n)} \mathbf{t}_2, n \in \mathbb{N},$

(**3.8**) $\mathbf{t}_1 \in_{w,[0]} \mathbf{t}_2, \mathbf{t}_1 \in_{w,[1]} \mathbf{t}_2, \ldots, \mathbf{t}_1 \in_{w,[n]} \mathbf{t}_2, n \in \mathbb{N},$

where $\mathbf{t}_1$ and $\mathbf{t}_2$ are numerical terms and $X$ is any set variable.

*The intended meanings of these respective atomic formulas are that:*

(**3.1**) $\mathbf{t}_1$ equals $\mathbf{t}_2$ in a *strong consistent sense*,
$\mathbf{t}_1$ is less than $\mathbf{t}_2$, in a *strong consistent sense*,
$\mathbf{t}_1$ equals $\mathbf{t}_2$ in a *weak inconsistent sense*,
$\mathbf{t}_1$ is less than $\mathbf{t}_2$, in a *weak inconsistent sense*;

(**3.2**) $\mathbf{t}_1$ equals $\mathbf{t}_2$ in a *weak inconsistent sense with rank* $= 0$,
$\mathbf{t}_1$ equals $\mathbf{t}_2$ in a *weak inconsistent sense with rank* $= 1,2,\ldots,$
$\mathbf{t}_1$ equals $\mathbf{t}_2$ in a *weak consistent sense with rank* $= n, n \in \mathbb{N}$;

(**3.3**) $\mathbf{t}_1$ equals $\mathbf{t}_2$ in a *strictly inconsistent sense with rank* $= 0$,
$\mathbf{t}_1$ equals $\mathbf{t}_2$ in a *strictly inconsistent sense with rank* $= 1,\ldots,$
$\mathbf{t}_1$ equals $\mathbf{t}_2$ in a *strictly inconsistent sense with rank* $= n, n \in \mathbb{N}$;

(**3.4**) $\mathbf{t}_1$ is less than $\mathbf{t}_2$, in a *weak inconsistent sense with rank* $= 0$,

$\mathbf{t}_1$is less than $\mathbf{t}_2$, in a *weak inconsistent sense with rank* $= 1,2,\ldots,$

$\mathbf{t}_1$ is less than $\mathbf{t}_2$, in a *weak consistent sense with rank* $= n, n \in \mathbb{N}$;

(**3.5**) $\mathbf{t}_1$ is less than $\mathbf{t}_2$, in a *strictly inconsistent sense with rank* $= 0$,

$\mathbf{t}_1$is less than $\mathbf{t}_2$, in a *strictly inconsistent sense with rank* $= 1,2,\ldots,$

$\mathbf{t}_1$ is less than $\mathbf{t}_2$, in a *strictly inconsistent sense with rank* $= n, n \in \mathbb{N}$;

(**3.6**) $\mathbf{t}_1$ is an element of $X$ in a *strong consistent sense*,

$\mathbf{t}_1$ is an element of $X$ in a *weak inconsistent sense*,

(**3.7**) $\mathbf{t}_1$ is an element of $X$ in a *weak inconsistent sense with rank* $= 0$,

$\mathbf{t}_1$ is an element of $X$ in a *weak inconsistent sense* with *rank* $= 1,2\ldots,$

$\mathbf{t}_1$ is an element of $X$ in a *weak inconsistent sense* with *rank* $= n$;

(**3.8**) $\mathbf{t}_1$ is an element of $X$ in a *strictly inconsistent sense with rank* $= 0$,

$\mathbf{t}_1$ is an element of $X$ in a *strictly inconsistent sense* with *rank* $= 1,2\ldots,$

$\mathbf{t}_1$ is an element of $X$ in a *strictly inconsistent sense* with *rank* $= n, n \in \mathbb{N}$;

(**4**) Formulas are built up from:

(**4.1**) atomic formulas by means of propositional connectives $\wedge, \vee, \neg, \rightarrow, \leftrightarrow$ (and, or, not, implies, if and only if),

(**4.2**) *consistent and inconsistent number* quantifiers $\forall n, \exists n$ (for all $n$, there exists $n$),

(**4.3**) *consistent and inconsistent set* quantifiers $\forall X, \exists X$ (for all $X$, there exists $X$),

(**4.4**) operators $(\cdot)^{(0)}, (\cdot)^{(1)}, \ldots, (\cdot)^{(n)}, (\cdot)^{[0]}, (\cdot)^{[1]}, \ldots, (\cdot)^{[n]}, n \in \mathbb{N}$.

(**5**) A sentence is a formula with no free variables.

**Definition 3.1.3**. (language of second order inconsistent arithmetic). $L_2^{\#}$ is defined to be the language of second order inconsistent arithmetic as described above. In writing terms and formulas of $L_2^{\#}$, we shall use parentheses and brackets to indicate grouping, as is customary in mathematical logic textbooks.

The semantics of the paraconsistent language $L_2^{\#}$ are given by the following definition.

**Definition 3.1.4**. (paraconsistent $L_2^{\#}$-structures). A paraconsistent model for $L_2^{\#}$, also called a paraconsistent structure for $L_2^{\#}$ or an paraconsistent $L_2^{\#}$-structure, is an ordered $19$-tuple:

$$\mathbf{M}_{\mathbf{inc}} = \breve{\mathbf{M}} = \Big\{|\breve{\mathbf{M}}|, \mathbf{S}_{\breve{\mathbf{M}}}, \big(\cdot +_{\breve{\mathbf{M}}} \cdot\big), (\cdot \times_{\breve{\mathbf{M}}} \cdot), 0_{\mathbf{s}}^{\breve{\mathbf{M}}}, 1_{\mathbf{s}}^{\breve{\mathbf{M}}}, 0_{w}^{\breve{\mathbf{M}}}, 1_{w}^{\breve{\mathbf{M}}}, \Big\{0_{w,(n)}^{\breve{\mathbf{M}}}\Big\}_{n\in\mathbb{N}},$$

$$\Big\{0_{w,[n]}^{\breve{\mathbf{M}}}\Big\}_{n\in\mathbb{N}}, \Big\{1_{w,(n)}^{\breve{\mathbf{M}}}\Big\}_{n\in\mathbb{N}}, \Big\{1_{w,[n]}^{\breve{\mathbf{M}}}\Big\}_{n\in\mathbb{N}}, \big(\cdot =_{\mathbf{s}}^{\breve{\mathbf{M}}} \cdot\big), \big(\cdot =_{w}^{\breve{\mathbf{M}}} \cdot\big),$$

$$\Big\{\big(\cdot =_{w,(n)}^{\breve{\mathbf{M}}} \cdot\big)\Big\}_{n\in\mathbb{N}}, \Big\{\big(\cdot =_{w,[n]}^{\breve{\mathbf{M}}} \cdot\big)\Big\}_{n\in\mathbb{N}}, \tag{3.1.3}$$

$$\big(\cdot <_{\mathbf{s}}^{\breve{\mathbf{M}}} \cdot\big), \big(\cdot <_{w}^{\breve{\mathbf{M}}} \cdot\big), \Big\{\big(\cdot <_{w,(n)}^{\breve{\mathbf{M}}} \cdot\big)\Big\}_{n\in\mathbb{N}}, \Big\{\big(\cdot <_{w,[n]}^{\breve{\mathbf{M}}} \cdot\big)\Big\}_{n\in\mathbb{N}}\Big\},$$

$$\mathbb{N} = \{0, 1, 2, \ldots\},$$

where $|\mathbf{M}|_{\mathbf{inc}} = |\breve{\mathbf{M}}|$ is an *inconsistent set* which serves as the range of the consistent and inconsistent number variables,$\mathbf{S}_{\breve{\mathbf{M}}}$ is a set of subsets of $|\breve{\mathbf{M}}|$ serving as the range of the set variables, $+_{\breve{\mathbf{M}}}$ and $\times_{\breve{\mathbf{M}}}$ are binary operations on $|\breve{\mathbf{M}}|$, $0_{\breve{\mathbf{M}}} \triangleq \Big\{0_{\mathbf{s}}^{\breve{\mathbf{M}}}, 0_{w}^{\breve{\mathbf{M}}}, 0_{w,(n)}^{\breve{\mathbf{M}}}, 0_{w,[n]}^{\breve{\mathbf{M}}}\Big\}$ and $1_{\breve{\mathbf{M}}} \triangleq \Big\{1_{\mathbf{s}}^{\breve{\mathbf{M}}}, 1_{w}^{\breve{\mathbf{M}}}, 1_{w,(n)}^{\breve{\mathbf{M}}}, 1_{w,[n]}^{\breve{\mathbf{M}}}\Big\}, n \in \mathbb{N}$ are distinguished elements of $|\breve{\mathbf{M}}|$, and $\big(\cdot =_{\mathbf{s}}^{\breve{\mathbf{M}}} \cdot\big), \big(\cdot <_{\mathbf{s}}^{\breve{\mathbf{M}}} \cdot\big)$ is a binary *strongly consistent relations* on $|\breve{\mathbf{M}}|, \big(\cdot =_{w}^{\breve{\mathbf{M}}} \cdot\big), \big(\cdot =_{w,(n)}^{\breve{\mathbf{M}}} \cdot\big)\big(\cdot <_{w}^{\breve{\mathbf{M}}} \cdot\big), \big(\cdot <_{w,(n)}^{\breve{\mathbf{M}}} \cdot\big), n \in \mathbb{N}$ is a binary *weakly inconsistent relations* on $|\breve{\mathbf{M}}|$, $\big(\cdot =_{w,[n]}^{\breve{\mathbf{M}}} \cdot\big), \big(\cdot <_{w,[n]}^{\breve{\mathbf{M}}} \cdot\big), n \in \mathbb{N}$ is a binary *inconsistent relations* on $|\breve{\mathbf{M}}|$. We always assume that the *inconsistent* sets $|\breve{\mathbf{M}}|$ and $\mathbf{S}_{\breve{\mathbf{M}}}$ are disjoint and nonempty. Formulas of $L_2^{\#}$ are interpreted in *inconsistent set* $\breve{\mathbf{M}}$ in the obvious way.

**Definition 3.1.5**. Strictly $\in$-consistent ($\in_{\mathbf{s}}$ -consistent) set $X$ it is a set such that $\forall x(x \in_{\mathbf{s}} X \vee x \notin_{\mathbf{s}} X)$.

**Definition 3.1.6**. Weakly $\in$-inconsistent ($\in_{w}$ -inconsistent) set $X$ it is a set such that $\forall x(x \in_{w} X \vee x \notin_{w} X)$.

**Definition 3.1.7**. Weakly $\in$-inconsistent with rank$=n, n \in \mathbb{N}$ ($\in_{w,(n)}$ -inconsistent) set $X$ it is a set such that:

$$\forall x(x \in_{w,(n)} X \vee x \notin_{w,(n)} X). \tag{3.1.4}$$

**Definition 3.1.8.**Strictly $\in$-inconsistent with rank$=n, n \in \mathbb{N}$ ($\in_{w,[n]}$ -inconsistent) set $X$ it is a set such that:

$$\forall x(x \in_{w,[n]} X \wedge x \notin_{w,[n]} X). \tag{3.1.5}$$

**Definition 3.1.9**. An strictly $\in$-consistent set $\breve{\mathbf{A}} \subseteq_{\mathbf{s}} |\mathbf{M}|_{\mathbf{inc}} =_{\mathbf{s}} |\breve{\mathbf{M}}|$ is said to be strictly consistent definable over $\breve{\mathbf{M}}$ allowing parameters from $\breve{\mathbf{B}}$ if there exists a formula $\varphi(n)$ with parameters from $\breve{\mathbf{B}}$ and no free variables other than $n$ such that

$$\breve{\mathbf{A}} = \{a \in_{\mathbf{s}} |\breve{\mathbf{M}}| : \breve{\mathbf{M}} \models \varphi(a), (\varphi(a))^{[0]} \vdash \mathbf{C}\}. \tag{3.1.6}$$

Here $\breve{\mathbf{M}} \models \varphi(a)$ as it is usual means that $\breve{\mathbf{M}}$ satisfies $\varphi(a)$,i.e., $\varphi(a)$ is true in $\breve{\mathbf{M}}$ and $\mathbf{C}$ is an any sentence of the inconsistent language $L_2^{\#}$.

**Definition 3.1.10**.The *intended model* for $L_2^{\#}$ is the model

$$\breve{\mathbf{M}} \triangleq$$

$$\triangleq (\mathbb{N}^{\#}, \mathbf{P}(\mathbb{N}^{\#}), +, \times, 0_{\breve{\mathbf{M}}} \triangleq \{0_{\mathbf{s}}, 0_w, 0_{w,(i)}, 0_{w,[i]}\}, 1_{\breve{\mathbf{M}}} \triangleq \{1_{\mathbf{s}}, 1_w, 1_{w,(i)}, 1_{w,[i]}\}, \tag{3.1.6'}$$

$$<_{\breve{\mathbf{M}}} \triangleq \{<_{\mathbf{s}}, <_w, <_{w,(i)}, <_{w,[i]}\}, i \in \mathbb{N})$$

where $\mathbb{N}^{\#} \triangleq \mathbb{N}_{\mathbf{inc}}$ is the set of paraconsistent natural numbers, $\mathbf{P}(\mathbb{N}^{\#})$ is the set of all **s**-subsets of $\mathbb{N}^{\#}$, and $+, \times, 0_{\breve{\mathbf{M}}} \triangleq \{0_{\mathbf{s}}, 0_w, 0_{w,(i)}, 0_{w,[i]}\}, 1_{\breve{\mathbf{M}}} \triangleq \{1_{\mathbf{s}}, 1_w, 1_{w,(i)}, 1_{w,[i]}\}$, $<_{\breve{\mathbf{M}}} \triangleq \{<_{\mathbf{s}}, <_w, <_{w,(i)}, <_{w,[i]}\}, i \in \mathbb{N}$ are as below.

**Definition 3.1.11**. (second order paraconsistent arithmetic $\mathbf{Z}_2^{\#}$). The axioms of second order paraconsistent arithmetic $\mathbf{Z}_2^{\#}$ consist of the universal closures of the following $L_2^{\#}$-formulas:

(i) **basic axioms**:

(i.1.) **basic axioms of the first group**:

(**i.1.1.**) $n + 1_{\mathbf{s}} \neq_{\mathbf{s}} 0_{\mathbf{s}}$,

(**i.1.2.**) $n + 1_w \neq_{\mathbf{s}} 0_w$,

(**i.1.3.**) $n + 1_{w,(i)} \neq_{\mathbf{s}} 0_{w,(i)}, i \in \mathbb{N}$,

(**i.1.4.**) $n + 1_{w,[i]} \neq_{\mathbf{s}} 0_{w,[i]}, i \in \mathbb{N}$,

(**i.1.5.**) $1_{\mathbf{s}} =_{\mathbf{s}} 1_{\mathbf{s}}, 0_{\mathbf{s}} =_{\mathbf{s}} 0_{\mathbf{s}}$,

**(i.1.6.)** $1_w =_{\mathbf{s}} 1_w, 0_w =_{\mathbf{s}} 0_w$,

**(i.1.7.)** $1_{w,(i)} =_{\mathbf{s}} 1_{w,(i)}, 0_{w,(i)} =_{\mathbf{s}} 0_{w,(i)}, i \in \mathbb{N}$,

**(i.1.8.)** $1_{w,[i]} =_{\mathbf{s}} 1_{w,[i]}, 0_{w,[i]} =_{\mathbf{s}} 0_{w,[i]}, i \in \mathbb{N}$;

(**i.2.**) **basic axioms of the second group**:

**(i.2.1.)** $m + 1_{\mathbf{s}} =_{\mathbf{s}} n + 1_{\mathbf{s}} \to m =_{\mathbf{s}} n$,

**(i.2.2.)** $m + 1_w =_w n + 1_w \to m =_w n$,

**(i.2.3.)** $m + 1_{w,(i)} =_{w,(i)} n + 1_{w,(i)} \to m =_{w,(i)} n, i \in \mathbb{N}$,

**(i.2.4.)** $m + 1_{w,[i]} =_{w,[i]} n + 1_{w,[i]} \to m =_{w,[i]} n, i \in \mathbb{N}$;

(**i.3.**) **basic axioms of the third group**:

**(i.3.1.)** $m + 0_{\mathbf{s}} =_{\mathbf{s}} m$,

**(i.3.2.)** $m + 0_w =_w m$,

**(i.3.3.)** $m + 0_{w,(i)} =_{w,(i)} m, i \in \mathbb{N}$,

**(i.3.4.)** $m + 0_{w,[i]} =_{w,[i]} m, i \in \mathbb{N}$;

(**i.4.**) **basic axioms of the fourth group**:

**(i.4.1.)** $m + (n + 1_{\mathbf{s}}) =_{\mathbf{s}} (m + n) + 1_{\mathbf{s}}$,

**(i.4.2.)** $m + (n + 1_w) =_w (m + n) + 1_w$,

**(i.4.3.)** $m + (n + 1_{w,(i)}) =_{w,(i)} (m + n) + 1_{w,(i)}, i \in \mathbb{N}$,

**(i.4.4.)** $m + (n + 1_{w,[i]}) =_{w,[i]} (m + n) +_{w,[i]} 1_{w,[i]}, i \in \mathbb{N}$;

(**i.5.**) **basic axioms of the fifth group**:

**(i.5.1.)** $m \times 0_{\mathbf{s}} =_{\mathbf{s}} 0_{\mathbf{s}}$,

**(i.5.2.)** $m \times 0_w =_w 0_w$,

**(i.5.3.)** $m \times 0_{w,(i)} =_{w,(n)} 0_{w,(i)}, i \in \mathbb{N}$,

**(i.5.4.)** $m \times 0_{w,[i]} =_{w,[i]} 0_{w,[i]}, i \in \mathbb{N}$;

(**i.6.**) **basic axioms of the sixth group**:

**(i.6.1.)** $m \times (n + 1_{\mathbf{s}}) =_{\mathbf{s}} (m \times n) + m$,

**(i.6.2.)** $m \times (n + 1_w) =_w (m \times n) + m$,

**(i.6.3.)** $m \times (n + 1_{w,(i)}) =_{w,(i)} (m \times n) + m, i \in \mathbb{N}$,

**(i.6.4.)** $m \times (n + 1_{w,[i]}) =_{w,[i]} (m \times n) + m, i \in \mathbb{N}$;

(**i.7**) **basic axioms of the seventh group**:

**(i.7.1.)** $\neg(m <_{\mathbf{s}} 0_{\mathbf{s}})$,

**(i.7.2.)** $\neg(m <_w 0_w)$,

**(i.7.3.)** $\neg(m <_{w,(i)} 0_{w,(i)}), i \in \mathbb{N}$,

**(i.7.4.)** $\neg(m <_{w,[i]} 0_{w,[i]}), i \in \mathbb{N}$;

(**i.8.**) **basic axioms of the eighth group**:

**(i.8.1.)** $m <_{\mathbf{s}} n + 1_{\mathbf{s}} \leftrightarrow (m <_{\mathbf{s}} n \vee m =_{\mathbf{s}} n)$,

**(i.8.2.)** $m <_w n + 1_w \leftrightarrow (m <_w n \vee m =_w n)$,

(**i.8.3.**) $m <_{w,(i)} n+1_{w,(i)} \leftrightarrow (m <_{w,(i)} n \vee m =_{w,(i)} n), i \in \mathbb{N}$,
(**i.8.4.**) $m <_{w,[i]} n+1_{w,[i]} \leftrightarrow (m <_{w,[i]} n \vee m =_{w,[i]} n), i \in \mathbb{N}$;
**Notation.3.1.1.** $\forall n_{n \in_{\boldsymbol{\alpha}} Z}(n \in_{\boldsymbol{\alpha}} X) \Leftrightarrow \forall n[(n \in_{\boldsymbol{\alpha}} Z) \wedge (n \in_{\boldsymbol{\alpha}} X)]$.

(**ii**) **induction axioms**:
(**ii.1.**) **strictly consistent induction (s-induction) axiom**
$\exists Y_1 \forall X[((0_{\mathbf{s}} \in_{\mathbf{s}} X) \wedge \forall n_{n \in_{\mathbf{s}} Y_1}[(n \in_{\mathbf{s}} X \to n+1_{\mathbf{s}} \in_{\mathbf{s}} X)] \to \forall n_{n \in_{\mathbf{s}} Y_1}(n \in_{\mathbf{s}} X)]$,
(**ii.2.**) **weakly inconsistent (*w*-inconsistent) induction axiom (*w*-induction axiom)**:
$\exists Y_2 \forall X[(0_w \in_w X \wedge \forall n_{n \in_w Y_2}[(n \in_w X \to n+1_w \in_w X)] \to \forall n_{n \in_w Y_2}(n \in_w X)]$,
(**ii.3.**) **weakly inconsistent with rang**$=n, n \in \mathbb{N}$ ($\{w,(n)\}$**-inconsistent**) **induction axiom** ($\{w,(n)\}$**-induction axiom**):
$\forall i_{(i \in \mathbb{N})} \exists Y_3^i \forall X\Big[(0_{w,(i)} \in_{w,(i)} X \wedge \forall n_{n \in_{w,(i)} Y_3^i}[(n \in_{w,(i)} X \to n+1_{w,(i)} \in_{w,(i)} X)]$
$\to \forall n_{n \in_{w,(i)} Y_3^i}(n \in_{w,(i)} X)\Big]$,
(**ii.4.**) **strictly inconsistent with rang**$=n, n \in \mathbb{N}$ ($\{w,[n]\}$**-inconsistent**) **induction axiom** ($\{w,[n]\}$**-induction axiom**):
$\forall i_{(i \in \mathbb{N})} \exists Y_4^i \forall X\Big[(0_{w,[i]} \in_{w,[i]} X \wedge \forall n_{n \in_{w,[i]} Y_4^i}[(n \in_{w,[i]} X \to n+1_{w,[i]} \in_{w,[i]} X)]$
$\to \forall n_{n \in_{w,[i]} Y_4^i}(n \in_{w,[i]} X)\Big]$,
(**ii.5.**) $\{w,(\mathbb{N})\}$**-induction axiom**:
$\exists Y_5 \forall X[\forall i_{(i \in \mathbb{N})}[(0_{w,(i)} \in_{w,(i)} X \wedge \forall n_{n \in_{w,(i)} Y_5}[(n \in_{w,(i)} X \to n+1_{w,(i)} \in_{w,(i)} X)]]$
$\to \forall j_{(j \in \mathbb{N})} \forall n_{n \in_{w,(j)} Y_5}(n \in_{w,(j)} X)]$,
(**ii.6.**) $\{w,[\mathbb{N}]\}$**-induction axiom**:
$\exists Y_6 \forall X[\forall i_{(i \in \mathbb{N})}[(0_{w,[i]} \in_{w,[i]} X \wedge \forall n_{n \in_{w,[i]} Y_5}[(n \in_{w,[i]} X \to n+1_{w,[i]} \in_{w,[i]} X)]]$
$\to \forall j_{(j \in \mathbb{N})} \forall n_{n \in_{w,[j]} Y_5}(n \in_{w,[j]} X)]$,
(**ii.7.**) **global** para**consistent induction axiom**:
$\exists Y_* \forall X[((0_{\mathbf{s}} \in_{\mathbf{s}} X) \wedge (0_w \in_{\mathbf{s}} X) \wedge (\forall i_{(i \in \mathbb{N})}(0_{w,(i)} \in_{\mathbf{s}} X)) \wedge (\forall i_{(i \in \mathbb{N})}(0_{w,[i]} \in_{\mathbf{s}} X)) \wedge$
$\wedge \forall n_{n \in_{\mathbf{s}} Y_*}(n \in_{\mathbf{s}} X \to n+1_{\mathbf{s}} \in_{\mathbf{s}} X) \wedge \forall n_{n \in_{\mathbf{s}} Y_*}(n \in_{\mathbf{s}} X \to n+1_w \in_{\mathbf{s}} X) \wedge$
$\wedge \{\forall i_{(i \in \mathbb{N})} \forall n_{n \in_{\mathbf{s}} Y_*}(n \in_{\mathbf{s}} X \to n+1_{w,(i)} \in_{\mathbf{s}} X)\} \wedge \{\forall i_{(i \in \mathbb{N})} \forall n_{n \in_{\mathbf{s}} Y_*}(n \in_{\mathbf{s}} X \to n+1_{w,[i]} \in_{\mathbf{s}} X)\}$
$\to \forall n_{n \in_{\mathbf{s}} Y_*}(n \in_{\mathbf{s}} X)]$.
**Definition 3.1.12.**

$$Y_1 \triangleq \mathbb{N}_{\mathbf{s}}^{\#} = \mathbb{N},$$

$$Y_2 \triangleq \mathbb{N}_w^{\#},$$

$$Y_3^i \triangleq \mathbb{N}_{w,(i)}^{\#},$$

$$Y_4^i \triangleq \mathbb{N}_{w,[i]}^{\#}, \tag{3.1.7}$$

$$Y_5 \triangleq \mathbb{N}_{w,(\mathbb{N})}^{\#},$$

$$Y_6 \triangleq \mathbb{N}_{w,[\mathbb{N}]}^{\#},$$

$$Y_* \triangleq \mathbb{N}^{\#} \triangleq \mathbb{N}_{\mathbf{pc}}.$$

(**iii**) **paraconsistent order axioms**:

(**iii**.**1**.) **weakly inconsistent order axiom** ($w$-**order axiom**):
every nonempty $w$-subset $X \subseteq_w \mathbb{N}_w^{\#}$ has $w$-least element,i.e.
a least element relative to weakly inconsistent order $(\bullet <_w \bullet)$.
(**iii**.**2**.) **weakly inconsistent with *rank*** $=n, n \in \mathbb{N}$ **order axiom**
($\{w,(n)\}$-**order axiom**):
every nonempty $\{w,(n)\}$-subset $X \subseteq_{w,(n)} \mathbb{N}_{w,(n)}^{\#}$ has $\{w,(n)\}$-least element,i.e.
a least element relative to weakly inconsistent order $(\bullet <_{w,(n)} \bullet)$.
(**iii**.**3**.) **strictly inconsistent with *rank*** $= n, n \in \mathbb{N}$ **order axiom**
($\{w,[n]\}$-**order axiom**):
every nonempty $\{w,[n]\}$-subset $X \subseteq_{w,[n]} \mathbb{N}_{w,(n)}^{\#}$ has $\{w,[n]\}$-least element,i.e.
a least element relative to strictly inconsistent order $(\bullet <_{w,[n]} \bullet)$.
(**iv**) **restricted comprehension scheme**:
**strictly consistent comprehension** (**s-comprehension**) **scheme**:
(**iv**.**1**.) $\exists X \forall n_{n \in_{\mathbf{s}} \mathbb{N}_{\mathbf{s}}^{\#}}(n \in_{\mathbf{s}} X \leftrightarrow \varphi(n))$,
where $\varphi(n)$ is any formula of $L_2^{\#}$ in which $X$ does not occur freely.
(**v**.) **non-restricted comprehension schemes**:
(**v**.**1**.) **weakly inconsistent comprehension scheme**
($w$-**comprehension scheme**):
$\exists X \forall n_{n \in_w \mathbb{N}_w^{\#}}(n \in_w X \leftrightarrow \varphi(n,X))$,

(**v.2.**) **weakly inconsistent with *rank*** $=n, n \in \mathbb{N}$ **comprehension scheme** ($\{w,(n)\}$**-comprehension scheme**):
$\exists X \forall n(n \in_{w,(i)} X \leftrightarrow (\varphi(n,X))^{(i)}), i \in \mathbb{N},$
(**v.3.**) **strictly inconsistent with *rank*** $=n, n \in \mathbb{N}$ **comprehension scheme** ($\{w,[n]\}$**-comprehension scheme**) :
$\exists X \forall n(n \in_{w,[i]} X \leftrightarrow (\varphi(n,X))^{[i]}), i \in \mathbb{N},$
where $\varphi(n)$ is any formula of $L_2^{\#}$ .

# III.2.Paraconsistent Mathematics Within $\mathbf{Z}_2^{\#}$.

We now outline the development of certain portions of paraconsistent mathematics within $\mathbf{Z}_2^{\#}$.

**Definition 3.2.1.**If $X$ and $Y$ are set variables, we use:
(**i**) $X =_{\mathbf{s}} Y$ and $X \subseteq_{\mathbf{s}} Y$ as abbreviations for the formulas $\forall n(n \in_{\mathbf{s}} X \leftrightarrow n \in_{\mathbf{s}} Y)$ and $\forall n(n \in_{\mathbf{s}} X \rightarrow n \in_{\mathbf{s}} Y)$ respectively.
(**ii**) $X =_w Y$ and $X \subseteq_w Y$ as abbreviations for the formulas $\forall n(n \in_w X \leftrightarrow n \in_w Y)$ and $\forall n(n \in_w X \rightarrow n \in_w Y)$ respectively.
(**iii**) $X =_{w,(i)} Y$ and $X \subseteq_{w,(i)} Y, i \in \mathbb{N}$ as abbreviations for the formulas: $\forall n(n \in_{w,(i)} X \leftrightarrow n \in_{w,(i)} Y)$ and $\forall n(n \in_{w,(i)} X \rightarrow n \in_{w,(i)} Y)$ respectively.
(**iv**) $X =_{w,[i]} Y$ and $X \subseteq_{w,[i]} Y, i \in \mathbb{N}$ as abbreviations for the formulas: $\forall n(n \in_{w,[i]} X \leftrightarrow n \in_{w,[i]} Y)$ and $\forall n(n \in_{w,[i]} X \rightarrow n \in_{w,[i]} Y)$ respectively.
(**v**) a strictly $\in$-consistent ($\in_{\mathbf{s}}$ -consistent) set $X \subset_{\mathbf{s}} \mathbb{N}^{\#}$ is defined to be a set such that: $\forall n[(n \in_{\mathbf{s}} X) \vee (n \notin_{\mathbf{s}} X)]$
(**vi**) a weakly $\in$-inconsistent ($\in_w$ -inconsistent) set $X \subset_w \mathbb{N}^{\#}$ is defined to be a set such that: $\forall x(x \in_w X \vee x \notin_w X)$.

(**vii**) a weakly $\in$-inconsistent with rank $=n, n \in \mathbb{N}$ ($\in_{w,(n)}$ -inconsistent) set $X$ it is a set $X \subset_{w.(n)} \mathbb{N}^{\#}$ such that:
$$\forall x(x \in_{w,(n)} X \vee x \notin_{w,(n)} X)$$
(**viii**) a strictly $\in$-inconsistent with rank $=n, n \in \mathbb{N}$ ($\in_{w,[n]}$ -inconsistent) set $X$ it is a set $X \subset_{w,[n]} \mathbb{N}^{\#}$ such that:
$$\forall x((x \in_{w,[n]} X) \wedge (x \notin_{w,[n]} X))$$
**Definition 3.2.2**. Strictly consistent single element set:

$$\{x\}_{\mathbf{s}} \triangleq \forall y[y \in_{\mathbf{s}} \{x\}_{\mathbf{s}} \leftrightarrow y =_{\mathbf{s}} x]. \quad (3.2.1)$$

Weakly inconsistent single element set:

$$\{x\}_w \triangleq \forall y[y \in_w \{x\}_w \leftrightarrow y =_w x]. \qquad (3.2.2)$$

Weakly inconsistent with rank$=n, n \in \mathbb{N}$ single element set:

$$\{x\}_{w,(n)} \triangleq \forall y\left[y \in_{w,(n)} \{x\}_{w,(n)} \leftrightarrow y =_{w,(n)} x\right]. \qquad (3.2.3)$$

Strictly inconsistent with rank$= n, n \in \mathbb{N}$ single element set:

$$\{x\}_{w,[n]} \triangleq \forall y\left[y \in_{w,[n]} \{x\}_{w,[0]} \leftrightarrow y =_{w,[0]} x\right]. \qquad (3.2.4)$$

**Definition 3.2.3**. Strictly consistent two-element set:

$$\{x,y\}_{\mathbf{s}} \triangleq \forall z[z \in_{\mathbf{s}} \{x,y\}_{\mathbf{s}} \leftrightarrow (z =_{\mathbf{s}} x) \vee (z =_{\mathbf{s}} y)]. \qquad (3.2.5)$$

Weakly inconsistent two-element set:

$$\{x,y\}_w \triangleq \forall z[z \in_w \{x,y\}_w \leftrightarrow (z =_w x) \vee (z =_w y)]. \qquad (3.2.6)$$

Weakly inconsistent with rank$=n, n \in \mathbb{N}$ two-element set:

$$\{x,y\}_{w,(n)} \triangleq \forall z\left[z \in_{w,(n)} \{x,y\}_{w,(n)} \leftrightarrow (z =_{w,(n)} x) \vee (z =_{w,(n)} y)\right]. \qquad (3.2.6')$$

Strictly inconsistent with rank$=n, n \in \mathbb{N}$ two-element set:

$$\{x,y\}_{[n]} \triangleq \{x,y\}_{w,[n]} \triangleq \forall z\left[z \in_{w,[n]} \{x,y\}_{w,[n]} \leftrightarrow (z =_{w,[n]} x) \vee (z =_{w,[n]} y)\right]. \qquad (3.2.7)$$

**Definition 3.2.4**. Strictly consistent ordered pair:

$$(x,y)_{\mathbf{s}} \triangleq \left\{\{x\}_{\mathbf{s}},\ \{x,y\}_{\mathbf{s}}\right\}_{\mathbf{s}} \tag{3.2.8}$$

Weakly inconsistent ordered pair:

$$(x,y)_{w} \triangleq \left\{\{x\}_{w},\ \{x,y\}_{w}\right\}_{w} \tag{3.2.9}$$

Weakly inconsistent with rank$=n, n \in \mathbb{N}$ ordered pair:

$$(x,y)_{w,(n)} \triangleq \left\{\{x\}_{w,(n)},\ \{x,y\}_{w,(n)}\right\}_{w,(n)} \tag{3.2.10}$$

Strictly inconsistent with rank$=n, n \in \mathbb{N}$ ordered pair:

$$(x,y)_{w,[n]} \triangleq \left\{\{x\}_{w,[n]},\ \{x,y\}_{w,[n]}\right\}_{w,[n]} \tag{3.2.11}$$

**Definition 3.2.5**.The strictly consistent cartesian product of $X$ and $Y$ is:

$$X \times_{\mathbf{s}} Y \triangleq \{(x,y)_{\mathbf{s}} | (x \in_{\mathbf{s}} X) \wedge (y \in_{\mathbf{s}} Y)\}_{\mathbf{s}} \tag{3.2.12}$$

The weakly inconsistent cartesian product of $X$ and $Y$ is:

$$X \times_{w} Y \triangleq \{(x,y)_{w} | (x \in_{w} X) \wedge (y \in_{w} Y)\}_{w} \tag{3.2.13}$$

The weakly inconsistent with rank$=n, n \in \mathbb{N}$ cartesian product of $X$ and $Y$ is:

$$X \times_{w,(n)} Y \triangleq \left\{(x,y)_{w,(n)} | (x \in_{w,(n)} X) \wedge (y \in_{w,(n)} Y)\right\}_{w,(n)} \quad (3.2.14)$$

The strictly inconsistent cartesian product with rank$= n, n \in \mathbb{N}$ of $X$ and $Y$ is:

$$X \times_{w,[n]} Y \triangleq \left\{(x,y)_{w,[n]} | (x \in_{w,[n]} X)_{w,[n]} \wedge (y \in_{w,[n]} Y)_{w,[n]}\right\}_{w,[n]} \quad (3.2.15)$$

**Definition 3.2.6.**Within $\mathbf{Z}_2^{\#}$,we define $\mathbb{N}^{\#} \triangleq \mathbb{N}_{\mathbf{pc}}$ to be the unique set $X$ such that $\forall n_{n \in_{\mathbf{s}} \mathbb{N}^{\#}}(n \in_{\mathbf{s}} X)$.

**Definition 3.2.7.**(**i**) For $X, Y \subseteq_{\mathbf{s}} \mathbb{N}^{\#}$, a strictly consistent (**s**-consistent) function $f_{\mathbf{s}} : X \to_{\mathbf{s}} Y$ is defined to be a $\in$-consistent set $f \subseteq_{\mathbf{s}} X \times_{\mathbf{s}} Y$ such that for all $m \in_{\mathbf{s}} X$ there is exactly one $n \in_{\mathbf{s}} Y$ such that $(m,n)_{\mathbf{s}} \in_{\mathbf{s}} f$. For $m \in_{\mathbf{s}} X$, $f(m)$ is defined to be the $\mathbf{s}$-unique $n$ such that $(m,n)_{\mathbf{s}} \in_{\mathbf{s}} f$.
The usual properties of such functions can be proved in $\mathbf{Z}_2^{\#}$.
(**ii**) For $X, Y \subseteq_{w} \mathbb{N}^{\#}$, a weakly inconsistent ($w$-**inconsistent**) function $f_w : X \to_w Y$ is defined to be a weakly $\in$-inconsistent set $f_w \subseteq_w X \times_w Y$ such that for all $m$ such that $(m \in_{\mathbf{s}} X) \vee (m \in_w X)$ there is exist $w$-unique $n \in_w Y$ such that $(m,n)_w \in_w f_w$. For all $m$ such that $(m \in_{\mathbf{s}} X) \vee (m \in_w X)$, $f_w(m)$ is defined to be the $w$-unique $n$ such that $(m,n)_w \in_w f_w$.
(**iii**) For $X, Y \subseteq_{w,(n)} \mathbb{N}^{\#}_{w,(n)}$, a weakly inconsistent with rank$=n, n \in \mathbb{N}$ ($\{w,(n)\}$-**inconsistent**) function $f_{w,(n)} : X \to_{w,(n)} Y$ is defined to be a weakly (with rang$=n, n \in \mathbb{N}$) $\in_{w,(n)}$ - inconsistent set $f_{w,(n)} \subseteq_{w,(n)} X \times_{w,(n)} Y$ such that for all $m$ such that $(m \in_{\mathbf{s}} X) \vee (m \in_{w,(n)} X)$ there is exist $\{w,(n)\}$-unique $n \in_{w,(n)} Y$ such that $(m,n)_{w,(n)} \in_{w,(n)} f_{w,(n)}$. For $m \in_{w,(n)} X$, $f(m)$ is defined to be the $\{w,(n)\}$-unique $n$ such that $(m,n)_{w,(n)} \in_{w,(n)} f_{w,(n)}$.
(**iv**) For $X, Y \subseteq_{w,[n]} \mathbb{N}^{\#}$, a strictly inconsistent with rank$=n, n \in \mathbb{N}$ ($\{w,[n]\}$-**inconsistent**) function $f_{w,[n]} : X \to_{w,[n]} Y$ is defined to be a strictly inconsistent set $f_{w,[n]} \subseteq_{w,[n]} X \times_{w,[n]} Y$ such that for all $m \in_{w,[n]} X$ there is exist $\{w,[n]\}$-unique $n \in_{w,[n]} Y$ such that $(m,n)_{w,[n]} \in_{w,[n]} f$. For $m \in_{w,[n]} X$, $f(m)$ is defined to be the $\{w,[n]\}$-unique $n$ such that $(m,n)_{w,[n]} \in_{w,[n]} f_{w,[n]}$.

**Definition 3.2.8.**(**i**) (**the strictly consistent primitive recursion** (**s-recursion**)). This means that,given $f_{\mathbf{s}} : X \to_{\mathbf{s}} Y$ and $g_{\mathbf{s}} : \mathbb{N}^{\#}_{\mathbf{s}} \times_{\mathbf{s}} X \times_{\mathbf{s}} Y \to Y$, there is a $\mathbf{s}$-unique $h_{\mathbf{s}} : \mathbb{N}^{\#}_{\mathbf{s}} \times_{\mathbf{s}} X \to_{\mathbf{s}} Y$ defined by $h_{\mathbf{s}}(0_{\mathbf{s}}, m) = f_{\mathbf{s}}(m)$, $h_{\mathbf{s}}(n + 1_{\mathbf{s}}, m) = g_{\mathbf{s}}(n, m, h_{\mathbf{s}}(n, m))$ for all $n \in_{\mathbf{s}} \mathbb{N}^{\#}_{\mathbf{s}}$ and $m \in_{\mathbf{s}} X$.
The existence of $h_{\mathbf{s}}$ is proved by strictly consistent arithmetical comprehension, and the $\mathbf{s}$-uniqueness of $h_{\mathbf{s}}$ is proved by strictly consistent arithmetical induction.
(**ii**) (**the weakly consistent primitive recursion** ($w$-**recursion**)).This means that, given $f_w : X \to_w Y$ and $g_w : \mathbb{N} \times_w X \times_w Y \to_w Y,$ there is a unique $h_{\mathbf{s}} : \mathbb{N}^{\#}_w \times_w X \to_w Y$ defined by $h_w(0_w,m) = f_w(m), h(n+1_w,m) = g_w(n,m,h_w(n,m))$ for all $n \in_w \mathbb{N}^{\#}_w$ and $m \in_w X$.
The existence of $h_w$ is proved by weakly consistent arithmetical comprehension, and the $w$-uniqueness of $h_w$ is proved by weakly consistent arithmetical induction.
(**iii**) (**the weakly consistent with rank**$=i, i \in \mathbb{N}$ **primitive recursion** ($\{w,(n)\}$-**recursion**) ).This means that, given $f_{w,(n)} : X \to_{w,(n)} Y$ and $g_{w,(n)} : \mathbb{N}^{\#}_{w,(n)} \times_{w,(n)} X \times_{w,(n)} Y \to_{w,(n)} Y,$ there is a unique $h_{\mathbf{s}} : \mathbb{N}^{\#}_{w,(n)} \times_{w,(n)} X \to_{w,(n)} Y$ defined by $h_{w,(n)}(0_{w,(n)},m) = f_{w,(n)}(m), h(n+1_{w,(n)},m) = g_{w,(n)}(n,m,h_{w,(n)}(n,m))$ for all $n \in_{w,(n)} \mathbb{N}^{\#}_{w,(n)}$ and $m \in_{w,(n)} X$.
The existence of $h_{w,(n)}$ is proved by weakly consistent with rank$=n, n \in \mathbb{N}$ arithmetical comprehension and the uniqueness of $h_{w,(n)}$ is proved by weakly consistent with rank$=n, n \in \mathbb{N}$ arithmetical induction.
(**iv**) (**the strictly inconsistent with rank**$=i, i \in \omega$ **primitive recursion** ($\{w,[n]\}$-**recursion**)).This means that, given $f_{w,[n]} : X \to_{w,[n]} Y$ and $g_{w,[n]} : \mathbb{N} \times_{w,[n]} X \times_{w,[n]} Y \to_{w,[n]} Y,$ there is a unique $h_{\mathbf{s}} : \mathbb{N}^{\#}_{w,[n]} \times_{w,[n]} X \to_{w,[n]} Y$ defined by $h_{w,[n]}(0_{w,[n]},m) =_{w,[n]} f_{w,[n]}(m), h(n+1_{w,[n]},m) =_{w,[n]} g_{w,[n]}(n,m,h_{w,[n]}(n,m))$ for all $n \in_{w,[n]} \mathbb{N}^{\#}_{w,[n]}$ and $m \in_{w,[n]} X$.
The existence of $h_{w,[n]}$ is proved by strictly inconsistent with rank$=n, n \in \mathbb{N}$ arithmetical comprehension and the uniqueness of $h_{w,[n]}$ is proved by strictly inconsistent with rank$=n, n \in \omega$ arithmetical induction.
(**v**) (**the global paraconsistent primitive recursion**).This means that, given $f^{\mathbf{s}}_{\mathbf{gl}} : X \to_{\mathbf{s}} Y,\ f^{w}_{\mathbf{gl}} : X \to_w Y,\ f^{w.(i)}_{\mathbf{gl}} : X \to_{\mathbf{s}} Y,\ f^{w,[i]}_{\mathbf{gl}} : X \to_{\mathbf{s}} Y, i \in \mathbb{N}$ and given $g^{\mathbf{s}}_{\mathbf{gl}} : \mathbb{N}^{\#} \times_{\mathbf{s}} X \times_{\mathbf{s}} Y \to_{\mathbf{s}} Y, g^{w}_{\mathbf{gl}} : \mathbb{N}^{\#} \times_{\mathbf{s}} X \times_{\mathbf{s}} Y \to_{\mathbf{s}} Y,$ $g^{w,(i)}_{\mathbf{gl}} : \mathbb{N}^{\#} \times_{\mathbf{s}} X \times_{\mathbf{s}} Y \to_{\mathbf{s}} Y, g^{w,[i]}_{\mathbf{gl}} : \mathbb{N}^{\#} \times_{\mathbf{s}} X \times_{\mathbf{s}} Y \to_{\mathbf{s}} Y, i \in \mathbb{N}$ there is a weakly unique $h_{\mathbf{gl}} : \mathbb{N}^{\#} \times_{\mathbf{s}} X \to_{\mathbf{s}} Y$ defined by:
$h^{\mathbf{s}}_{\mathbf{gl}}(0_{\mathbf{s}},m) =_{\mathbf{s}} f^{\mathbf{s}}_{\mathbf{s}}(m), h^{w}_{\mathbf{gl}}(0_w,m) =_{\mathbf{s}} f^{w}_{\mathbf{gl}}(m),$
$h^{w,(i)}_{\mathbf{gl}}(0_{w,(i)},m) =_{\mathbf{s}} f^{w,(i)}_{\mathbf{gl}}(m),$
$h_{\mathbf{gl}}(0_{w,[i]},m) =_{\mathbf{s}} f^{w,[i]}_{\mathbf{gl}}(m), i \in \mathbb{N},$
$h^{\mathbf{s}}_{\mathbf{gl}}(n+1_{\mathbf{s}},m) =_{\mathbf{s}} g^{\mathbf{s}}_{\mathbf{gl}}(n,m,h_{\mathbf{s}}(n,m)),$
$h_{\mathbf{gl}}(n+1_w,m) =_{\mathbf{s}} g^{w}_{\mathbf{gl}}(n,m,h_{\mathbf{gl}}(n,m)),$
$h_{\mathbf{gl}}(n+1_{w,(i)},m) =_{\mathbf{s}} g^{w,(i)}_{\mathbf{gl}}(n,m,h_{\mathbf{gl}}(n,m)),$
$h_{\mathbf{gl}}(n+1_{w,[n]},m) =_{\mathbf{s}} g^{w,[i]}_{\mathbf{gl}}(n,m,h_{\mathbf{s}}(n,m)),$
$i \in \mathbb{N},$ for all $n \in_{\mathbf{s}} \mathbb{N}^{\#}$ and $m \in_{\mathbf{s}} X$.
The existence of $h_{\mathbf{s}}$ is proved by strictly consistent arithmetical comprehension and the uniqueness of $h_{\mathbf{s}}$ is proved by global paraconsistent arithmetical induction.

**Remark.3.1.1.**In particular, we have:

(**i**) the strictly consistent exponential function $\exp(m,n)_{\mathbf{s}} =_{\mathbf{s}} (m)^n_{\mathbf{s}}$, defined by $(m)^{0_{\mathbf{s}}}_{\mathbf{s}} =_{\mathbf{s}} 1_{\mathbf{s}}$,
$(m)^{n+1_{\mathbf{s}}}_{\mathbf{s}} =_{\mathbf{s}} (m)^n_{\mathbf{s}} \times_{\mathbf{s}} m$ for all $m,n \in_{\mathbf{s}} \mathbb{N}^{\#}_{\mathbf{s}}$.

(**ii**) the weakly inconsistent exponential function $\exp_w(m,n) =_w (m)^n_w$, defined by $(m)^{0_w}_w =_{\mathbf{s}} 1_w, (m)^{n+1_w}_w =_w (m)^n_w \times_w m$ for all $m,n \in_{\mathbf{s}} \mathbb{N}^{\#}_w$.

(**iii**) the weakly inconsistent with rank$=i, i \in \mathbb{N}$ exponential function $\exp_{w,(i)}(m,n) =_w (m)^n_{w,(i)}$,defined by $(m)^{0_{w,(i)}}_{w,(i)} =_{\mathbf{s}} 1_{w,(i)}$,
$(m)^{n+1_{w,(i)}}_{w,(i)} =_{w,(i)} (m)^n_{w,(i)} \times_w m$ for all $m,n \in_{w,(i)} \mathbb{N}^{\#}_{w,(i)}$.

(**iv**) the strictly inconsistent with rank$=i, i \in \mathbb{N}$ exponential function $\exp_{w,[i]}(m,n) =_{w,[i]} (m)^n_{w,[i]}$, defined by $(m)^{0_{w,[i]}}_{w,[i]} =_{w,[i]} 1_{w,[i]}$,
$(m)^{n+1_{w,[i]}}_{w,[i]} =_{w,[i]} (m)^n_{w,[i]} \times m$ for all $m,n \in_{w,[i]} \mathbb{N}^{\#}_{w,[i]}$.

The usual properties of the exponential function can be proved in $\mathbf{Z}^{\#}_2$.

(**v**) the global paraconsistent exponential function $\exp_{\mathbf{gl}}(m,n) \triangleq (m)^n_{\mathbf{gl}}$, defined by: $(m)^{0_{\mathbf{s}}}_{\mathbf{gl}} =_{\mathbf{s}} 1_{\mathbf{s}}, (m)^{n+1_{\mathbf{s}}}_{\mathbf{gl}} =_{\mathbf{s}} (m)^n_{\mathbf{gl}} \times m$ for all $m,n \in_{\mathbf{s}} \mathbb{N}^{\#}_{\mathbf{s}}$,
$(m)^{0_w}_{\mathbf{gl}} =_w 1_w, (m)^{n+1_w}_{\mathbf{gl}} =_w (m)^n_{\mathbf{gl}} \times m$ for all $m,n \in_{\mathbf{s}} \mathbb{N}^{\#}_w$,
$(m)^{0_{w,(n)}}_{\mathbf{gl}} =_{w,(n)} 1_{w,(n)}, (m)^{n+1_{w,(n)}}_{\mathbf{gl}} =_{w,(n)} (m)^n_{\mathbf{gl}} \times m$ for all $m,n \in_{\mathbf{s}} \mathbb{N}^{\#}_{w,(n)}$,
$(m)^{0_{w,[n]}}_{\mathbf{gl}} =_{w,[n]} 1_{w,[n]}, (m)^{n+1_{w,[n]}}_{\mathbf{gl}} =_{w,[n]} (m)^n_{\mathbf{gl}} \times m$ for all $m,n \in_{\mathbf{s}} \mathbb{N}^{\#}_{w,[n]}$.

The usual properties of the exponential function can be proved in $\mathbf{Z}^{\#}_2$.

Within $\mathbf{Z}^{\#}_2$, we define a numerical strictly consistent pairing function $\pi_{\mathbf{s}}(m,n)$ by $\pi_{\mathbf{s}}(m,n) =_{\mathbf{s}} (m+n)^2 + m$. Within $\mathbf{Z}^{\#}_2$ we can prove that, for all $m,n,i,j \in_{\mathbf{s}} \mathbb{N}^{\#}, \pi_{\mathbf{s}}(m,n) =_{\mathbf{s}} \pi_{\mathbf{s}}(i,j)$ if and only if $m =_{\mathbf{s}} i$ and $n =_{\mathbf{s}} j$. Moreover, using strictly consistent arithmetical comprehension,we can prove that for all sets $X, Y \subseteq_{\mathbf{s}} \mathbb{N}^{\#}$,there exists a set $\pi_{\mathbf{s}}(X \times Y) \subseteq_{\mathbf{s}} \mathbb{N}^{\#}$ consisting of all $\pi_{\mathbf{s}}(m,n)$ such that $m \in_{\mathbf{s}} X$ and $n \in_{\mathbf{s}} Y$. In particular we have $\pi_{\mathbf{s}}(\mathbb{N}^{\#} \times \mathbb{N}^{\#}) \subseteq_{\mathbf{s}} \mathbb{N}^{\#}$.

The paraconsistent natural number system is essentially already given to us by the language $L^{\#}_2$ and axioms of $\mathbf{Z}^{\#}_2$.Thus, within $\mathbf{Z}^{\#}_2$,a consistent and inconsistent natural number is defined to be an element of $\mathbb{N}^{\#}$,and the *paraconsistent natural number system* is defined to be the paraconsistent structure:

$$\mathbb{N}^{\#}, +_{\mathbb{N}^{\#}}, \times_{\mathbb{N}^{\#}}, 0_{\mathbb{N}^{\#}}, 1_{\mathbb{N}^{\#}}, <_{\mathbb{N}^{\#}}, =_{\mathbb{N}^{\#}},$$

$$0_{\mathbb{N}^{\#}} \triangleq \{0_{\mathbf{s}}, 0_{w}, 0_{w,(i)}, 0_{w,[i]}\},$$

$$1_{\mathbb{N}^{\#}} \triangleq \{1_{\mathbf{s}}, 1_{w}, 1_{w,(i)}, 1_{w,[i]}\}, \quad (3.2.16)$$

$$<_{\mathbb{N}^{\#}} \triangleq\triangleq \{<_{\mathbf{s}}, <_{w}, <_{w,(i)}, <_{w,[i]}\}$$

$$=_{\mathbb{N}^{\#}} \triangleq\triangleq \{=_{\mathbf{s}}, =_{w}, =_{w,(i)}, =_{w,[i]}\}$$

where $+_{\mathbb{N}^{\#}} : \mathbb{N}^{\#} \times_{\mathbf{s}} \mathbb{N}^{\#} \to_{\mathbf{s}} \mathbb{N}^{\#}$ is defined by $m +_{\mathbb{N}^{\#}} n =_{\mathbf{s}} m + n$, etc.
Thus for instance $+_{\mathbb{N}^{\#}}$ is the set of triples $((m,n)_{\mathbf{s}}, k)_{\mathbf{s}} \in_{\mathbf{s}} (\mathbb{N}^{\#} \times_{\mathbf{s}} \mathbb{N}^{\#}) \times_{\mathbf{s}} \mathbb{N}^{\#}$ such that $m + n =_{\mathbf{s}} k$. The existence of this set follows from the strictly consistent arithmetical comprehension.
In a standard manner, we can define within $\mathbf{Z}_2^{\#}$ the set $\mathbb{Z}^{\#} \triangleq \mathbb{Z}_{\mathbf{pc}}$ of the all consistent and inconsistent integers and the set of the all consistent and inconsistent rational numbers: $\mathbb{Q}^{\#} \triangleq \mathbb{Q}_{\mathbf{pc}}$.

**Definition 3.2.9.**(paraconsistent rational numbers $\mathbb{Q}^{\#}$) Let

$$(\mathbb{Z}^{\#})^{+} = \{a \in_{\mathbf{s}} \mathbb{Z}^{\#} : 0_{\mathbf{s}} <_{\mathbb{Z}^{\#}} a\}_{\mathbf{s}} \quad (3.2.17)$$

be the set of positive consistent integers,and let
(**a**) $\equiv_{\mathbb{Q}^{\#},\mathbf{s}}$ be the strictly equivalence relation on $\mathbb{Z}^{\#} \times_{\mathbf{s}} (\mathbb{Z}^{\#})^{+}$
defined by $(a,b)_{\mathbf{s}} \equiv_{\mathbb{Q}^{\#}} (c,d)_{\mathbf{s}}$ if and only if $a \times_{\mathbb{Z}^{\#}} d =_{\mathbf{s}} b \times_{\mathbb{Z}^{\#}} c$.
Then $\mathbb{Q}^{\#}$ is defined to be the set of all $(a,b)_{\mathbf{s}} \in \mathbb{Z}^{\#} \times_{\mathbf{s}} (\mathbb{Z}^{\#})^{+}$ such
that $(a,b)_{\mathbf{s}}$ is the $<_{\mathbf{s}}^{\mathbb{N}^{\#}}$ -minimum element of its $\equiv_{\mathbb{Q}^{\#}}$ -equivalence class.
Operations $+_{\mathbb{Q}^{\#}}, -_{\mathbb{Q}^{\#}}, \times_{\mathbb{Q}^{\#}}$ on $\mathbb{Q}^{\#}$ are defined by:
$(a,b)_{\mathbf{s}} +_{\mathbb{Q}^{\#}} (c,d)_{\mathbf{s}} \equiv_{\mathbb{Q}^{\#}} (a \times_{\mathbb{Z}^{\#}} d +_{\mathbb{Z}^{\#}} b \times_{\mathbb{Z}^{\#}} c, b \times_{\mathbb{Z}^{\#}} d)_{\mathbf{s}}$,
$-_{\mathbb{Q}^{\#}} (a,b)_{\mathbf{s}} \equiv_{\mathbb{Q}^{\#}} (-_{\mathbb{Z}^{\#}} a, b)_{\mathbf{s}}$, and
$(a,b)_{\mathbf{s}} \times_{\mathbb{Q}^{\#}} (c,d)_{\mathbf{s}} \equiv Q(a \times_{\mathbb{Z}^{\#}} c, b \times_{\mathbb{Z}^{\#}} d)_{\mathbf{s}}$. We let $0_{\mathbb{Q}^{\#}} \equiv_{\mathbb{Q}^{\#}} (0_{\mathbb{Z}^{\#}}, 1_{\mathbb{Z}^{\#}})$, i.e.
$0_{\mathbf{s},\mathbb{Q}^{\#}} \equiv_{\mathbb{Q}^{\#}} (0_{\mathbf{s},\mathbb{Z}^{\#}}, 1_{\mathbf{s},,\mathbb{Z}^{\#}}), 0_{w,\mathbb{Q}^{\#}} \equiv_{\mathbb{Q}^{\#}} (0_{w,\mathbb{Z}^{\#}}, 1_{w,\mathbb{Z}^{\#}}), 0_{w,(n),\mathbb{Q}^{\#}} \equiv_{\mathbb{Q}^{\#}} (0_{w,(n),\mathbb{Z}^{\#}}, 1_{w,(n),\mathbb{Z}^{\#}})$,
$0_{w,[n],\mathbb{Q}^{\#}} \equiv_{\mathbb{Q}^{\#}} (0_{w,[n],\mathbb{Z}^{\#}}, 1_{w,[n],\mathbb{Z}^{\#}})$,
and
$1_{\mathbb{Q}^{\#}} \equiv_{\mathbb{Q}^{\#}} (1_{\mathbb{Z}^{\#}}, 1_{\mathbb{Z}^{\#}})_{\mathbf{s}}$, i.e. $1_{\mathbf{s},\mathbb{Q}^{\#}} \equiv_{\mathbb{Q}^{\#}} (1_{\mathbf{s},\mathbb{Z}^{\#}}, 1_{\mathbf{s},\mathbb{Z}^{\#}})_{\mathbf{s}}, 1_{w,\mathbb{Q}^{\#}} \equiv_{\mathbb{Q}^{\#}} (1_{w,\mathbb{Z}^{\#}}, 1_{w,\mathbb{Z}^{\#}})_{\mathbf{s}}$,

$1_{w,(n),\mathbb{Q}^\#} \equiv_{\mathbb{Q}^\#} (1_{w,(n),\mathbb{Z}^\#}, 1_{w,(n),\mathbb{Z}^\#})_{\mathbf{s}}, 1_{w,[n],\mathbb{Q}^\#} \equiv_{\mathbb{Q}^\#} (1_{w,[n],\mathbb{Z}^\#}, 1_{w,[n],\mathbb{Z}^\#})_{\mathbf{s}},$
and we define a binary relations $<_{\mathbb{Q}^\#} \triangleq \{<_{\mathbf{s}}, <_w, <_{w,(i)}, <_{w,[i]}\}$ on $\mathbb{Q}^\#$ by letting $(a,b)_{\mathbf{s}} <_{\mathbb{Q}^\#} (c,d)_{\mathbf{s}}$ if
and only if $a \times_{\mathbb{Z}^\#} d <_{\mathbb{Z}^\#} b \times_{\mathbb{Z}^\#} c$. Finally $=_{\mathbb{Q}^\#} \triangleq \{=_{\mathbf{s}}, =_w, =_{w,(i)}, =_{w,[i]}\}$ is the identity relations on $\mathbb{Q}^\#$. We can then prove within $\mathbf{Z}_2^\#$ that the paraconsistent rational number system

$$\{\mathbb{Q}^\#, +_{\mathbb{Q}^\#}, -_{\mathbb{Q}^\#}, \times_{\mathbb{Q}^\#}, 0_{\mathbb{Q}^\#}, 1_{\mathbb{Q}^\#}, <_{\mathbb{Q}^\#}, =_{\mathbb{Q}^\#}\},$$

$$0_{\mathbb{Q}^\#} \triangleq \{0_{\mathbf{s},\mathbb{Q}^\#}, 0_{w,\mathbb{Q}^\#}, 0_{w,(n),\mathbb{Q}^\#}, 0_{w,[n],\mathbb{Q}^\#}\},$$

$$1_{\mathbb{Q}^\#} \triangleq \{1_{\mathbf{s},\mathbb{Q}^\#}, 1_{w,\mathbb{Q}^\#}, 1_{w,(n),\mathbb{Q}^\#}, 1_{w.[n],\mathbb{Q}^\#}\}, \quad (3.2.18)$$

$$<_{\mathbb{Q}^\#} \triangleq \{<_{\mathbf{s},\mathbb{Q}^\#}, <_{w,\mathbb{Q}^\#}, <_{w,(i),\mathbb{Q}^\#}, <_{w,[i],\mathbb{Q}^\#}\},$$

$$=_{\mathbb{Q}^\#} \triangleq \{=_{\mathbf{s},\mathbb{Q}^\#}, =_{w.\mathbb{Q}^\#}, =_{w,(i),\mathbb{Q}^\#}, =_{w,[i],\mathbb{Q}^\#}\},$$

has the usual properties of an paraordered field, etc.
We make the usual identifications whereby $\mathbb{N}^\#$ is regarded as a subset of $\mathbb{Z}^\#$ and $\mathbb{Z}^\#$ is regarded as a subset of $\mathbb{Q}^\#$. Namely $m \in_{\mathbf{s}} \mathbb{N}^\#$ is identified with $(m,0)_{\mathbf{s}} \in_{\mathbf{s}} \mathbb{Z}^\#$, and $a \in_{\mathbf{s}} \mathbb{Z}^\#$ is identified with $(a, 1_{\mathbb{Z}^\#})_{\mathbf{s}} \in_{\mathbf{s}} \mathbb{Q}^\#$.
We use $+$ ambiguously to denote $+_{\mathbb{N}^\#}, +_{\mathbb{Z}^\#}$, or $+_{\mathbb{Q}^\#}$ and similarly for $-, \times, 0, 1, <$. For $q, r \in_{\mathbf{s}} \mathbb{Q}^\#$ we write $q - r = q + (-r)$, and if $r \neq_{\mathbb{Q}^\#} 0$, $q/r =_{\mathbb{Q}^\#}$ the unique $q' \in_{\mathbf{s}} \mathbb{Q}^\#$ such that $q = q' \times r$. The function $\mathbf{s}\text{-}\exp(q,a) = \mathbf{s}\text{-}q^a$ for $q \in_{\mathbf{s}} \mathbb{Q}^\# \backslash \{0_{\mathbb{Q}^\#}\}$ and $a \in_{\mathbf{s}} \mathbb{Z}^\#$ is obtained by primitive recursion in the obvious way.

**Definition 3.2.10**.The absolute value functions:
(**i**) the strictly consistent absolute value function:
$|\cdot|_{\mathbf{s}} : \mathbb{Q}^\# \to_{\mathbf{s}} \mathbb{Q}^\#$ is defined by $|q|_{\mathbf{s}} =_{\mathbf{s}} q$ if $0 \leq_{\mathbf{s}} q$, and by $|q|_{\mathbf{s}} =_{\mathbf{s}} -q$ otherwise,
(**ii**) the weakly inconsistent absolute value function:
$|\cdot|_w : \mathbb{Q}^\# \to_{\mathbf{s}} \mathbb{Q}^\#$ is defined by $|q|_w =_{\mathbf{s}} q$ if $0 \leq_w q$, and by $|q|_w =_{\mathbf{s}} -q$ otherwise,
(**iii**) the weakly inconsistent with rang$= n, n \in \omega$ absolute value function:
$|\cdot|_{w,(n)} : \mathbb{Q}^\# \to_{\mathbf{s}} \mathbb{Q}^\#$ is defined by $|q|_{w,(n)} =_{\mathbf{s}} q$ if $0 \leq_{w,(n)} q$, and by $|q|_w =_{\mathbf{s}} -q$ otherwise,
(**iv**) the strictly inconsistent with rang$= n, n \in \omega$ absolute value function:
$|\cdot|_{w,[n]} : \mathbb{Q}^\# \to_{\mathbf{s}} \mathbb{Q}^\#$ is defined by $|q|_{w,[n]} =_{\mathbf{s}} q$ if $0 \leq_{w,[n]} q$, and by $|q|_{w,[n]} =_{\mathbf{s}} -q$ otherwise.

**Definition 3.2.11**.(**a**) A strictly consistent sequence of paraconsistent rational numbers is defined to be a strictly consistent function $f_{\mathbf{s}} : \mathbb{N}^{\#} \to_{\mathbf{s}} \mathbb{Q}^{\#}$.
We denote such a strictly consistent sequence as $\langle q_n : n \in_{\mathbf{s}} \mathbb{N}^{\#}\rangle_{\mathbf{s}}$, or simply $\langle q_n\rangle_{\mathbf{s}}$, where $q_n =_{\mathbf{s}} f_{\mathbf{s}}(n)$.

**Definition 3.2.12**. A double strictly consistent sequence of paraconsistent rational numbers is a consistent function $f_{\mathbf{s}} : \mathbb{N}^{\#} \times_{\mathbf{s}} \mathbb{N}^{\#} \to \mathbb{Q}^{\#}$, denoted $\langle q_{mn} : m,n \in_{\mathbf{s}} \mathbb{N}^{\#}\rangle_{\mathbf{s}}$ or simply $\langle q_{mn}\rangle_{\mathbf{s}}$, where $q_{mn} =_{\mathbf{s}} f_{\mathbf{s}}(m,n)$.

**Definition 3.2.13**.(paraconsistent real numbers).(**a**) Within $\mathbf{Z}_2^{\#}$, a strictly consistent real number is defined to be a Cauchy strictly consistent sequence of paraconsistent rational numbers,i.e.,as a strictly consistent sequence of paraconsistent rational numbers $x =_{\mathbf{s}} \langle q_n : n \in_{\mathbf{s}} \mathbb{N}^{\#}\rangle_{\mathbf{s}}$ such that:

$$\forall\varepsilon(\varepsilon \in_{\mathbf{s}} \mathbb{Q}^{\#})(0_{\mathbf{s}} <_{\mathbf{s}} \varepsilon \to \exists m\forall n(m <_{\mathbf{s}} n \to |q_m - q_n|_{\mathbf{s}} <_{\mathbf{s}} \varepsilon)). \qquad (3.2.19)$$

(**b**) Within $\mathbf{Z}_2^{\#}$, a weakly inconsistent real number is defined to be a Cauchy strictly consistent sequence of paraconsistent rational numbers, i.e.,as a strictly consistent sequence of paraconsistent rational numbers $x =_{\mathbf{s}} \langle q_n : n \in_{\mathbf{s}} \mathbb{N}^{\#}\rangle_{\mathbf{s}}$ such that:

$$\forall\varepsilon(\varepsilon \in_{\mathbf{s}} \mathbb{Q}^{\#})(0_w <_w \varepsilon \to \exists m\forall n(m <_w n \to |q_m - q_n|_w <_w \varepsilon)). \qquad (3.2.20)$$

(**c**) Within $\mathbf{Z}_2^{\#}$, a weakly inconsistent with rang$= n, n \in \omega$ real number is defined to be a Cauchy strictly consistent sequence of paraconsistent rational numbers,i.e.,as a strictly consistent sequence of paraconsistent rational numbers $x =_{\mathbf{s}} \langle q_n : n \in_{\mathbf{s}} \mathbb{N}^{\#}\rangle_{\mathbf{s}}$ such that:

$$\forall\varepsilon(\varepsilon \in_{\mathbf{s}} \mathbb{Q}^{\#})(0_{w,(n)} <_{w,(n)} \varepsilon \to \exists m\forall n(m <_{w,(n)} n \to |q_m - q_n|_{w,(n)} <_{w,(n)} \varepsilon)). \qquad (3.2.21)$$

(**d**) Within $\mathbf{Z}_2^{\#}$, a strictly inconsistent with rang$= n, n \in \omega$ real number is defined to be a Cauchy strictly consistent sequence of paraconsistent rational numbers,i.e.,as a strictly consistent sequence of paraconsistent rational numbers $x =_{\mathbf{s}} \langle q_n : n \in_{\mathbf{s}} \mathbb{N}^{\#}\rangle_{\mathbf{s}}$ such that:

$$\forall\varepsilon(\varepsilon \in_{\mathbf{s}} \mathbb{Q}^{\#})(0_{w,[n]} <_{w,[n]} \varepsilon \to \exists m\forall n(m <_{w,[n]} n \to |q_m - q_n|_{w,[n]} <_{w,[n]} \varepsilon)). \quad (3.2.22)$$

**Definition 3.2.14.**(**a**) If $x =_{\mathbf{s}} q_n$ and $y =_{\mathbf{s}} q'_n$ are strictly consistent real numbers, we write $x =_{\mathbf{s},\mathbb{R}^{\#}} y$ to mean that $\mathbf{s}\text{-}\lim_n |q_n - q'_n|_{\mathbf{s}} =_{\mathbf{s}} 0_{\mathbf{s}}$, i.e.,

$$\forall\varepsilon(\varepsilon \in_{\mathbf{s}} \mathbb{Q}^{\#})(0_{\mathbf{s}} <_{\mathbf{s}} \varepsilon \to \exists m\forall n(m <_{\mathbf{s}} n \to |q_n - q'_n|_{\mathbf{s}} <_{\mathbf{s}} \varepsilon)), \quad (3.2.23)$$

and we write $x <_{\mathbf{s},\mathbb{R}^{\#}} y$ to mean that

$$\exists\varepsilon(0 <_{\mathbf{s},\mathbb{R}^{\#}} \varepsilon \wedge \exists m\forall n(m <_{\mathbf{s}} n \to q_n + \varepsilon <_{\mathbf{s}} q'_n)). \quad (3.2.24)$$

(**b**) If $x =_w q_n$ and $y =_w q'_n$ are weakly consistent real numbers, we write $x =_{w,\mathbb{R}^{\#}} y$ to mean that $w\text{-}\lim_n |q_n - q'_n|_w =_w 0_w$, i.e.,

$$\forall\varepsilon(\varepsilon \in_{\mathbf{s}} \mathbb{Q}^{\#})(0_w <_w \varepsilon \to \exists m\forall n(m <_{\mathbf{s}} n \to |q_n - q'_n|_w <_w \varepsilon)), \quad (3.2.25)$$

and we write $x <_{w,\mathbb{R}^{\#}} y$ to mean that

$$\exists\varepsilon(0 <_w \varepsilon \wedge \exists m\forall n(m <_{\mathbf{s}} n \to q_n + \varepsilon <_w q'_n)). \quad (3.2.26)$$

(**c**) If $x =_{w,(n)} q_n$ and $y =_{w,(n)} q'_n$ are weakly inconsistent with rang$= n, n \in \omega$ real numbers, we write $x =_{w,(n),\mathbb{R}^{\#}} y$ to mean that $\{w,(n)\}\text{-}\lim_n |q_n - q'_n|_{w,(n)} =_{w,(n)} 0_{w,(n)}$, i.e.,

$$\forall\varepsilon(\varepsilon \in_{\mathbf{s}} \mathbb{Q}^{\#})(0_{w,(n)} <_{w,(n)} \varepsilon \to \exists m\forall n(m <_{\mathbf{s}} n \to |q_n - q'_n|_{w,(n)} <_{w,(n)} \varepsilon)), \quad (3.2.27)$$

and we write $x <_{w,(n),\mathbb{R}^{\#}} y$ to mean that

$$\exists\varepsilon(0 <_{w,(n)} \varepsilon \wedge \exists m \forall n(m <_{\mathbf{s}} n \to q_n + \varepsilon <_{w,(n)} q_n')). \tag{3.2.28}$$

(**d**) If $x =_{\mathbf{s}} q_n$ and $y =_{\mathbf{s}} q_n'$ are strictly inconsistent with rang= $n, n \in \omega$ real numbers, we write $x =_{\mathbb{R}^\#,w,[n]} y$ to mean that $w,[n]\text{-}\lim_n |q_n - q_n'|_{w,[n]} =_{w,[n]} 0_{w,[n]}$, i.e.,

$$\forall\varepsilon(\varepsilon \in_{\mathbf{s}} \mathbb{Q}^\#)(0_{w,[n]} <_{w,[n]} \varepsilon \to \exists m \forall n(m <_{\mathbf{s}} n \to |q_n - q_n'|_{w,[n]} <_{w,[n]} \varepsilon)), \tag{3.2.29}$$

and we write $x <_{w,[n],\mathbb{R}^\#} y$ to mean that

$$\exists\varepsilon(0 <_{w,[n]} \varepsilon \wedge \exists m \forall n(m <_{\mathbf{s}} n \to q_n + \varepsilon <_{w,[n]} q_n')). \tag{3.2.30}$$

Also $x +_{\mathbb{R}^\#} y = \langle q_n + q_n' \rangle_{\mathbf{s}}, x \times_{\mathbb{R}^\#} y =_{\mathbf{s}} \langle q_n \times q_n' \rangle_{\mathbf{s}}, -_{\mathbb{R}^\#} x =_{\mathbf{s}} \langle -q_n \rangle_{\mathbf{s}}$, $0_{\mathbf{s},\mathbb{R}^\#} = \langle 0_{\mathbf{s}} \rangle_{\mathbf{s}}, 1_{\mathbf{s},\mathbb{R}^\#} = \langle 1_{\mathbf{s}} \rangle_{\mathbf{s}}$, etc.
We use $\mathbb{R}^\#$ or $\mathbb{R}^\#_{\mathbf{inc}}$ to denote the set of all *paraconsistent real numbers*.
Thus $x \in_{\mathbf{s}} \mathbb{R}^\#$ means that $x$ is *consistent or paraconsistent real number.* (Formally, we cannot speak of the set $\mathbb{R}^\#$ within the language of second order arithmetic, since it is a set of sets.)
We shall usually omit the subscript $\mathbb{R}^\#$ in $+_{\mathbb{R}^\#}, -_{\mathbb{R}^\#}, \times_{\mathbb{R}^\#}, 0_{\mathbf{s},\mathbb{R}^\#}, 0_{w,\mathbb{R}^\#}$, $0_{w,(n),\mathbb{R}^\#}, 0_{w,[n],\mathbb{R}^\#}, 1_{\mathbf{s},\mathbb{R}^\#}, 1_{w,\mathbb{R}^\#}, 1_{w,(n),\mathbb{R}^\#}, 1_{w,[n],\mathbb{R}^\#}, <_{\mathbf{s},\mathbb{R}^\#}, <_{w,\mathbb{R}^\#}\ <_{w,(n),\mathbb{R}^\#}, <_{w,[n],\mathbb{R}^\#}$, $=_{\mathbf{s},\mathbb{R}^\#}, =_{w,\mathbb{R}^\#}, =_{w,(n),\mathbb{R}^\#}, =_{w,[n],\mathbb{R}^\#}$.
Thus the *paraconsistent real number system* consists of $\mathbb{R}^\#, +, -, \times, 0_{\mathbf{s}}, 0_w$, $0_{w,(n)}, 0_{w,[n]}, 1_{\mathbf{s}}, 1_w, 1_{w,(n)}, 1_{w,[n]}, <_{\mathbf{s}}, <_w, <_{w,(n)}, <_{w,[n]}, =_{\mathbf{s}}, =_w, =_{w,(n)}, =_{w,[n]}$.
We shall sometimes identify a paraconsistent rational number $q \in_{\mathbf{s}} \mathbb{Q}^\#$ with the corresponding paraconsistent real number $x_q =_{\mathbf{s}} \langle q \rangle_{\mathbf{s}} \in_{\mathbf{s}} \mathbb{R}^\#$.
Within $\mathbf{Z}_2^\#$ one can prove that the real number system has the usual properties of an *paraconsistent Archimedean paraordered field*, etc. The complex paraconsistent numbers can be introduced as usual as pairs of paraconsistent real numbers.
Within $\mathbf{Z}_2^\#$,it is straightforward to carry out the proofs of all the basic results in real and complex *paraconsistent linear* and *paraconsistent polynomial* algebra. For example, the paraconsistent analog of the fundamental theorem of algebra can be proved in $\mathbf{Z}_2^\#$.
**Definition 3.2.15**. A strictly consistent sequence of paraconsistent real numbers is defined to be a double strictly consistent sequence of rational numbers $\langle q_{mn} : m,n \in_{\mathbf{s}} \mathbb{N}^{\#}\rangle_{\mathbf{s}}$ such that for each $m, \langle q_{mn} : n \in_{\mathbf{s}} \mathbb{N}^{\#}\rangle_{\mathbf{s}}$ is a consistent real number.

Such a strictly consistent sequence of paraconsistent real numbers is denoted $\langle x_m : m \in_{\mathbf{s}} \mathbb{N}^{\#}\rangle_{\mathbf{s}}$, where $x_m =_{\mathbf{s}} \langle q_{mn} : n \in_{\mathbf{s}} \mathbb{N}^{\#}\rangle$. Within $\mathbf{Z}_2^{\#}$ we can prove that every bounded (in *paraconsistent sense*) strictly consistent sequence of paraconsistent real numbers has a *paraconsistent least upper bound.* This is a very useful *paracompleteness* property of the paraconsistent real number system.

For instance, it implies that an infinite series of positive terms is *paraconvergent* if and only if the finite partial sums are *parabounded.*

We now turn of certain portions of paraconsistent abstract algebra within $\mathbf{Z}_2^{\#}$. Because of the restriction to the language $L_2^{\#}$ of second order paraconsistent arithmetic, we cannot expect to obtain a good general theory of arbitrary (countable and uncountable) paraconsistent algebraic structures. However, we can develop paraconsistent countable algebra, i.e., the theory of countable paraconsistent algebraic structures, within $\mathbf{Z}_2^{\#}$.

**Definition 3.2.16**. A countable paraconsistent commutative ring is defined within $\mathbf{Z}_2^{\#}$ to be a paraconsistent structure $\mathbf{R_{inc}}, +_{\mathbf{R_{inc}}}, -_{\mathbf{R_{inc}}}, \times_{\mathbf{R_{inc}}},$ $0_{\mathbf{s},\mathbf{R_{inc}}}, 0_{w,\mathbf{R_{inc}}}, 0_{w,(n),\mathbf{R_{inc}}}, 0_{w,[n],\mathbf{R_{inc}}}, 1_{\mathbf{s},\mathbf{R_{inc}}}, 1_{w,\mathbf{R_{inc}}}, 1_{w,(n),\mathbf{R_{inc}}}, 1_{w,[n],\mathbf{R_{inc}}}$ where $\mathbf{R_{inc}} \subseteq_{\mathbf{s}} \mathbb{N}^{\#}, +_{\mathbf{R_{inc}}} : \mathbf{R_{inc}} \times_{\mathbf{s}} \mathbf{R_{inc}} \to_{\mathbf{s}} \mathbf{R_{inc}},$ etc., and the usual commutative paraconsistent ring axioms are assumed.

(We include $0_{\mathbf{s}} \neq_{\mathbf{s}} 1_{\mathbf{s}}, 0_w \neq_{\mathbf{s}} 1_w, 0_{w,(n)} \neq_{\mathbf{s}} 1_{w,(n)}, 0_{w,[n]} \neq_{\mathbf{s}} 1_{w,[n]},$ among those axioms.) The subscript $\mathbf{R_{inc}}$ is usually omitted.

An strictly consistent ideal in $\mathbf{R_{inc}}$ is a set $I_{\mathbf{s}}^{\#} \subseteq_{\mathbf{s}} \mathbf{R_{inc}}$ such that $a \in_{\mathbf{s}} I_{\mathbf{s}}^{\#}$ and $b \in_{\mathbf{s}} I_{\mathbf{s}}^{\#}$ imply $a + b \in_{\mathbf{s}} I_{\mathbf{s}}^{\#}; a \in_{\mathbf{s}} I_{\mathbf{s}}^{\#}$ and $r \in_{\mathbf{s}} \mathbf{R_{inc}}$ imply $a \times r \in_{\mathbf{s}} I_{\mathbf{s}}^{\#}$, and $0_{\mathbf{s}} \in I_{\mathbf{s}}^{\#}$ and $1_{\mathbf{s}} \notin I_{\mathbf{s}}^{\#}$.

We define: (**a**) an strongly consistent equivalence relation $=_{I_{\mathbf{s}}^{\#},\mathbf{s}}$ on $\mathbf{R_{inc}}$ by $r =_{I_{\mathbf{s}}^{\#}} s$ if and only if $r - s \in_{\mathbf{s}} I_{\mathbf{s}}^{\#}$.

(**b**) an weakly inconsistent equivalence relation $=_{I_{\mathbf{s}}^{\#},w}$ on $\mathbf{R_{inc}}$ by $r =_{I_{\mathbf{s}}^{\#},w} s$ if and only if $r - s \in_w I_{\mathbf{s}}^{\#}$.

We let $\mathbf{R_{inc}}/I_{\mathbf{s}}^{\#}$ be the set of $r \in_{\mathbf{s}} \mathbf{R_{inc}}$ such that $r$ is the $<_{\mathbf{s},\mathbb{N}^{\#}}$ -minimum element of its equivalence class under $=_{I_{\mathbf{s}}^{\#},\mathbf{s}}$ .Thus $\mathbf{R_{inc}}/I_{\mathbf{s}}^{\#}$ consists of one element of each $=_{I_{\mathbf{s}}^{\#},\mathbf{s}}$ -equivalence class of elements of $\mathbf{R_{inc}}$. With the appropriate operations,$\mathbf{R}/I_{\mathbf{s}}^{\#}$ becomes a countable commutative ring, the quotient ring of $\mathbf{R_{inc}}$ by $I_{\mathbf{s}}^{\#}$. The ideal $I_{\mathbf{s}}^{\#}$ is said to be prime if $\mathbf{R_{inc}}/I_{\mathbf{s}}^{\#}$ is an integral domain, and maximal if $\mathbf{R_{inc}}/I_{\mathbf{s}}^{\#}$ is a field.

Next we indicate how some basic concepts and results of analysis and topology can be developed within $\mathbf{Z}_2^{\#}$.

**Definition 3.2.17**.Within $\mathbf{Z}_2^{\#}$,a **s**-paracomplete separable paraconsistent metric space is a nonempty set $A \subseteq_{\mathbf{s}} \mathbb{N}^{\#}$ together with a function $d_{\mathbf{s}} : A \times_{\mathbf{s}} A \to_{\mathbf{s}} \mathbb{R}^{\#}$ satisfying $a =_{\mathbf{s}} a \to d_{\mathbf{s}}(a,a) =_{\mathbf{s}} 0_{\mathbf{s}}$,
$0_{\mathbf{s}} \leq_{\mathbf{s}} d_{\mathbf{s}}(a,b) =_{\mathbf{s}} d_{\mathbf{s}}(b,a)$,
and $d_{\mathbf{s}}(a,c) \leq_{\mathbf{s}} d_{\mathbf{s}}(a,b) + d_{\mathbf{s}}(b,c)$ for all $a,b,c \in_{\mathbf{s}} A$.
(Formally, $d_{\mathbf{s}}$ is a strictly consistent sequence of paraconsistent real numbers, indexed by $A \times_{\mathbf{s}} A$.)
We define a point of the **s**-paracomplete separable paraconsistent metric space $\widehat{A}$ to be a sequence $x =_{\mathbf{s}} \langle a_n : n \in_{\mathbf{s}} \mathbb{N}^{\#} \rangle_{\mathbf{s}}$, $a_n \in_{\mathbf{s}} A$, satisfying

$$\forall \varepsilon(\varepsilon \in_{\mathbf{s}} \mathbb{R}^{\#})(0_{\mathbf{s}} <_{\mathbf{s}} \varepsilon \to \exists m \forall n(m <_{\mathbf{s}} n \to d_{\mathbf{s}}(a_m,a_n) <_{\mathbf{s}} \varepsilon)). \qquad (3.2.31)$$

The pseudometric $d_{\mathbf{s}}$ is extended from $A$ to $\left(\widehat{A}\right)_{\mathbf{s}}$ by

$$d_{\mathbf{s}}(x,y) = \mathbf{s}\text{-} \lim_{n\to\infty} d_{\mathbf{s}}(a_n,b_n) \qquad (3.2.32)$$

where $x =_{\mathbf{s}} \langle a_n : n \in \mathbb{N}^{\#} \rangle_{\mathbf{s}}$ and $y =_{\mathbf{s}} \langle b_n : n \in_{\mathbf{s}} \mathbb{N}^{\#} \rangle_{\mathbf{s}}$. We write $x =_{\mathbf{s}} y$ if and only if $d_{\mathbf{s}}(x,y) =_{\mathbf{s}} 0_{\mathbf{s}}$. For example, $\mathbb{R}^{\#} =_{\mathbf{s}} \widehat{\mathbb{Q}^{\#}}$ under the metric $d_{\mathbf{s}}(q,q') =_{\mathbf{s}} |q-q'|_{\mathbf{s}}$.

**Definition 3.2.18**.Within $\mathbf{Z}_2^{\#}$,a *weakly paracomplete* ($w$-paracomplete) separable paraconsistent metric space is a nonempty set $A \subseteq_{\mathbf{s}} \mathbb{N}^{\#}$ together with a function $d_w : A \times_{\mathbf{s}} A \to_{\mathbf{s}} \mathbb{R}^{\#}$ satisfying
$a =_{\mathbf{s}} a \to d_w(a,a) =_{\mathbf{s}} 0_{\mathbf{s}}, a =_w a \to d_w(a,a) =_w 0_w$,
$0_w \leq_w d_w(a,b) =_w d_w(b,a)$, and $d_w(a,c) \leq_w d_w(a,b) + d_w(b,c)$ for all $a,b,c \in_{\mathbf{s}} A$.
(Formally, $d_w$ is a consistent sequence of paraconsistent real numbers, indexed by $A \times_{\mathbf{s}} A$.)
We define a point of the $w$-paracomplete separable paraconsistent metric space $\left(\widehat{A}\right)_w$ to be a sequence $x =_{\mathbf{s}} \langle a_n : n \in_{\mathbf{s}} \mathbb{N}^{\#} \rangle_{\mathbf{s}}$, $a_n \in_{\mathbf{s}} A$, satisfying

$$\forall \varepsilon(\varepsilon \in_{\mathbf{s}} \mathbb{R}^{\#})(0_w <_w \varepsilon \to \exists m \forall n(m <_{\mathbf{s}} n \to d_w(a_m,a_n) <_w \varepsilon)). \qquad (3.2.33)$$

The pseudometric $d_w$ is extended from $A$ to $\left(\widehat{A}\right)_w$ by

$$d_w(x,y) = w\text{-}\lim_{n\to\infty} d_w(a_n,b_n) \tag{3.2.34}$$

where $x =_{\mathbf{s}} \langle a_n : n \in \mathbb{N}^{\#}\rangle_{\mathbf{s}}$ and $y =_{\mathbf{s}} \langle b_n : n \in_{\mathbf{s}} \mathbb{N}^{\#}\rangle_{\mathbf{s}}$. We write
(a) $x =_{\mathbf{s}} y$ if and only if $d_{\mathbf{s}}(x,y) =_{\mathbf{s}} 0_{\mathbf{s}}$,
(b) $x =_w y$ if and only if $d_w(x,y) =_w 0_w$.
For example, $\mathbb{R}^{\#}_w =_{\mathbf{s}} \left(\widehat{\mathbb{Q}^{\#}}\right)_w$ under the metric $d_w(q,q') =_{\mathbf{s}} |q-q'|_w$.

**Definition 3.2.19**.Within $\mathbf{Z}_2^{\#}$, a *weakly paracomplete* with rang$= n, n \in \omega$ ($\{w,(n)\}$-paracomplete) separable paraconsistent metric space is a nonempty set $A \subseteq_{\mathbf{s}} \mathbb{N}^{\#}$ together with a function $d_w : A \times_{\mathbf{s}} A \to_{\mathbf{s}} \mathbb{R}^{\#}$ satisfying:
$a =_{\mathbf{s}} a \to d_{w,(n)}(a,a) =_{\mathbf{s}} 0_{\mathbf{s}}, a =_{w,(n)} a \to d_{w,(n)}(a,a) =_{w,(n)} 0_{w,(n)}$,
$0_{w,(n)} \leq_{w,(n)} d_{w,(n)}(a,b) =_{w,(n)} d_{w,(n)}(b,a)$, and
$d_{w,(n)}(a,c) \leq_{w,(n)} d_{w,(n)}(a,b) + d_{w,(n)}(b,c)$
for all $a,b,c \in_{\mathbf{s}} A$.
(Formally, $d_{w,(n)}$ is a consistent sequence of paraconsistent real numbers, indexed by $A \times_{\mathbf{s}} A$.)
We define a point of the $\{w,(n)\}$-paracomplete separable paraconsistent metric space $\left(\widehat{A}\right)_{w,(n)}$ to be a sequence $x =_{\mathbf{s}} \langle a_n : n \in_{\mathbf{s}} \mathbb{N}^{\#}\rangle_{\mathbf{s}}$, $a_n \in_{\mathbf{s}} A$, satisfying:

$$\forall\varepsilon(\varepsilon \in_{\mathbf{s}} \mathbb{R}^{\#})(0_{w,(n)} <_{w,(n)} \varepsilon \to \exists m\forall n(m <_{\mathbf{s}} n \to d_{w,(n)}(a_m,a_n) <_{w,(n)} \varepsilon)). \tag{3.2.35}$$

**Definition 3.2.20**.Within $\mathbf{Z}_2^{\#}$, a *strictly paracomplete* with rang$= n, n \in \omega$ ($\{w,[n]\}$-paracomplete) separable paraconsistent metric space is a nonempty set $A \subseteq_{\mathbf{s}} \mathbb{N}^{\#}$ together with a function $d_w : A \times_{\mathbf{s}} A \to_{\mathbf{s}} \mathbb{R}^{\#}$ satisfying:
$a =_{\mathbf{s}} a \to d_{w,[n]}(a,a) =_{\mathbf{s}} 0_{\mathbf{s}}, a =_{w,[n]} a \to d_{w,[n]}(a,a) =_{w,[n]} 0_{w,[n]}$,
$0_{w,[n]} \leq_{w,[n]} d_{w,[n]}(a,b) =_{w,[n]} d_{w,[n]}(b,a)$, and
$d_{w,[n]}(a,c) \leq_{w,[n]} d_{w,[n]}(a,b) + d_{w,[n]}(b,c)$
for all $a,b,c \in_{\mathbf{s}} A$.
(Formally, $d_{w,[n]}$ is a consistent sequence of paraconsistent real numbers, indexed by $A \times_{\mathbf{s}} A$.)

We define a point of the $\{w,[n]\}$-paracomplete separable paraconsistent metric space $\left(\widehat{A}\right)_{w,[n]}$ to be a sequence $x =_{\mathbf{s}} \langle a_n : n \in_{\mathbf{s}} \mathbb{N}^{\#}\rangle_{\mathbf{s}}$, $a_n \in_{\mathbf{s}} A$, satisfying:

$$\forall\varepsilon(\varepsilon \in_{\mathbf{s}} \mathbb{R}^{\#})(0_{w,[n]} <_{w,[n]} \varepsilon \to \exists m\forall n(m <_{\mathbf{s}} n \to d_{w,[n]}(a_m,a_n) <_{w,[n]} \varepsilon)). \qquad (3.2.35)$$

**Definition 3.2.21**.(paraconsistent **s**-continuous functions). Within $\mathbf{Z}_2^{\#}$, if $\widehat{A}$ and $\widehat{B}$ are complete separable paraconsistent metric spaces,a paraconsistent **s**-continuous function $\phi : \widehat{A} \to_{\mathbf{s}} \widehat{B}$ is a set $\Phi_{\mathbf{s}} \subseteq_{\mathbf{s}} A \times_{\mathbf{s}} (\mathbb{Q}^{\#})^{+} \times_{\mathbf{s}} B \times_{\mathbf{s}} (\mathbb{Q}^{\#})^{+}$ satisfying the following coherence conditions:

$$1. [(a,r,b,s)_{\mathbf{s}} \in_{\mathbf{s}} \Phi_{\mathbf{s}}] \wedge [(a,r,b',s')_{\mathbf{s}} \in_{\mathbf{s}} \Phi_{\mathbf{s}}] \dashrightarrow d_{\mathbf{s}}(b,b') <_{\mathbf{s}} s+s';$$

$$2. [(a,r,b,s)_{\mathbf{s}} \in_{\mathbf{s}} \Phi_{\mathbf{s}}] \wedge [d_{\mathbf{s}}(b,b') + s <_{\mathbf{s}} s'] \dashrightarrow (a,r,b',s')_{\mathbf{s}} \in_{\mathbf{s}} \Phi_{\mathbf{s}} \qquad (3.2.36)$$

$$3. [(a,r,b,s)_{\mathbf{s}} \in_{\mathbf{s}} \Phi_{\mathbf{s}}] \wedge [d_{\mathbf{s}}(a,a') + r' <_{\mathbf{s}} r] \dashrightarrow (a',r',b,s)_{\mathbf{s}} \in_{\mathbf{s}} \Phi_{\mathbf{s}}$$

**Definition 3.2.22**.(paraconsistent $w$-paracontinuous functions). Within $\mathbf{Z}_2^{\#}$, if $\widehat{A}$ and $\widehat{B}$ are $w$-paracomplete separable paraconsistent metric spaces,a paraconsistent $w$-paracontinuous function $\phi : \widehat{A} \to_{\mathbf{s}} \widehat{B}$ is a set $\Phi_w \subseteq_w A \times_{\mathbf{s}} (\mathbb{Q}^{\#})^{+} \times_{\mathbf{s}} B \times_{\mathbf{s}} (\mathbb{Q}^{\#})^{+}$ satisfying the following coherence conditions:

$$1. [(a,r,b,s)_w \in_w \Phi_w] \wedge [(a,r,b',s')_w \in_w \Phi_w] \dashrightarrow d_w(b,b') <_w s+s';$$

$$2. [(a,r,b,s)_w \in_w \Phi_w] \wedge [d_w(b,b') + s <_w s'] \dashrightarrow (a,r,b',s')_w \in_w \Phi_w \qquad (3.2.37)$$

$$3. [(a,r,b,s)_w \in_w \Phi_w] \wedge [d_w(a,a') + r' <_w r] \dashrightarrow (a',r',b,s)_w \in_w \Phi_w$$

**Definition 3.2.23**.(paraconsistent $\{w,(n)\}$-paracontinuous functions).

Within $\mathbf{Z}_2^{\#}$,if $\widehat{A}$ and $\widehat{B}$ are $\{w,(n)\}$-paracomplete separable paraconsistent metric spaces,a paraconsistent $\{w,(n)\}$-paracontinuous function $\phi : \widehat{A} \to_{\mathbf{s}} \widehat{B}$ is a set $\Phi_w \subseteq_w A \times_{\mathbf{s}} (\mathbb{Q}^{\#})^+ \times_{\mathbf{s}} B \times_{\mathbf{s}} (\mathbb{Q}^{\#})^+$ satisfying the following coherence conditions:

$$
\begin{gathered}
1.\,[(a,r,b,s)_{w,(n)} \in_{w,(n)} \Phi_{w,(n)}] \wedge [(a,r,b',s')_{w,(n)} \in_{w,(n)} \Phi_{w,(n)}] \dashrightarrow \\
\to d_{w,(n)}(b,b') <_{w,(n)} s+s'; \\
2.\,[(a,r,b,s)_{w,(n)} \in_{w,(n)} \Phi_{w,(n)}] \wedge [d_{w,(n)}(b,b')+s <_{w,(n)} s'] \dashrightarrow \\
\to (a,r,b',s')_{w,(n)} \in_{w,(n)} \Phi_{w,(n)} \\
3.\,[(a,r,b,s)_{w,(n)} \in_{w,(n)} \Phi_{w,(n)}] \wedge [d_{w,(n)}(a,a')+r' <_{w,(n)} r] \dashrightarrow \\
\to (a',r',b,s)_{w,(n)} \in_{w,(n)} \Phi_{w,(n)}
\end{gathered}
\tag{3.2.38}
$$

**Definition 3.2.24.**(paraconsistent $\{w,[n]\}$-paracontinuous functions). Within $\mathbf{Z}_2^{\#}$,if $\widehat{A}$ and $\widehat{B}$ are $\{w,[n]\}$-paracomplete separable paraconsistent metric spaces,a paraconsistent $\{w,[n]\}$-paracontinuous function $\phi : \widehat{A} \to_{\mathbf{s}} \widehat{B}$ is a set $\Phi_w \subseteq_w A \times_{\mathbf{s}} (\mathbb{Q}^{\#})^+ \times_{\mathbf{s}} B \times_{\mathbf{s}} (\mathbb{Q}^{\#})^+$ satisfying the following coherence conditions:

$$1.\,[(a,r,b,s)_{w,[n]} \in_{w,[n]} \Phi_{w,[n]}] \wedge [(a,r,b',s')_{w,[n]} \in_{w,[n]} \Phi_{w,[n]}] \to$$

$$\to d_{w,[n]}(b,b') <_{w,[n]} s+s';$$

$$2.\,[(a,r,b,s)_{w,[n]} \in_{w,[n]} \Phi_{w,[n]}] \wedge [d_{w,[n]}(b,b')+s <_{w,[n]} s'] \to$$

$$\to (a,r,b',s')_{w,[n]} \in_{w,[n]} \Phi_{w,[n]} \quad (3.2.39)$$

$$3.\,[(a,r,b,s)_{w,[n]} \in_{w,[n]} \Phi_{w,[n]}] \wedge [d_{w,[n]}(a,a')+r' <_{w,[n]} r] \to$$

$$\to (a',r',b,s)_{w,[n]} \in_{w,[n]} \Phi_{w,[n]}$$

**Definition 3.2.25**. (paraconsistent **s**-open sets). Within $\mathbf{Z}_2^{\#}$, let $\widehat{A}$ be a **s**-paracomplete separable paraconsistent metric space. A (code for an) strictly open set (**s**-open set) in $\widehat{A}$ is any set $U \subseteq_{\mathbf{s}} \widehat{A} \times (\mathbb{Q}^{\#})^{+}$. For $x \in_{\mathbf{s}} \widehat{A}$ we write $x \in_{\mathbf{s}} U$ if and only if $d_{\mathbf{s}}(x,a) <_{\mathbf{s}} r$ for some $(a,r)_{\mathbf{s}} \in_{\mathbf{s}} U$.

**Definition 3.2.26**. (paraconsistent $w$-open sets). Within $\mathbf{Z}_2^{\#}$, let $\widehat{A}$ be a $w$-paracomplete separable paraconsistent metric space. A (code for an) weakly open set ($w$-open set) in $\widehat{A}$ is any set $U \subseteq_w \widehat{A} \times (\mathbb{Q}^{\#})^{+}$. For $x \in_w \widehat{A}$ we write $x \in_w U$ if and only if $d_w(x,a) <_w r$ for some $(a,r)_w \in_w U$.

**Definition 3.2.27**. (paraconsistent $\{w,(n)\}$-open sets). Within $\mathbf{Z}_2^{\#}$, let $\widehat{A}$ be a $\{w,(n)\}$-paracomplete separable paraconsistent metric space.

A (code for an) weakly open with rang$= n, n \in \omega$ set ($\{w,(n)\}$-open set) in $\widehat{A}$ is any set $U \subseteq_{w,(n)} \widehat{A} \times (\mathbb{Q}^{\#})^{+}$. For $x \in_{w,(n)} \widehat{A}$ we write $x \in_{w,(n)} U$ if and only if $d_{w,(n)}(x,a) <_{w,(n)} r$ for some $(a,r)_{w,(n)} \in_{w,(n)} U$.

**Definition 3.2.28**. (paraconsistent $\{w,[n]\}$-open sets). Within $\mathbf{Z}_2^{\#}$, let $\widehat{A}$ be a $\{w,[n]\}$-paracomplete separable paraconsistent metric space.

A (code for an) weakly open with rang$= n, n \in \omega$ set ($\{w,[n]\}$-open set) in $\widehat{A}$ is any set $U \subseteq_{w,[n]} \widehat{A} \times (\mathbb{Q}^{\#})^{+}$. For $x \in_{w,[n]} \widehat{A}$ we write $x \in_{w,[n]} U$ if and only if $d_{w,[n]}(x,a) <_{w,[n]} r$ for some $(a,r)_{w,[n]} \in_{w,[n]} U$.

**Definition 3.2.29**. A separable paraconsistent Banach **s**-space is defined within $\mathbf{Z}_2^{\#}$ to be a paracomplete separable metric space $\widehat{A}$ arising from a countable **s**-pseudonormed vector space $\widehat{A}$ over the paraconsistent field $\mathbb{Q}^{\#}$.
**Example 3.2.1**. (**a**) With the **s**-metric

$$d_{\mathbf{s}}(f,g) =_{\mathbf{s}} \mathbf{s}\text{-}\sup_{0_{\mathbf{s}} \leq_{\mathbf{s}} x \leq_{\mathbf{s}} 1_{\mathbf{s}}} |f_1(x) - f_2(x)|_{\mathbf{s}} \tag{3.2.40}$$

we have **s**-paracomplete separable paraconsistent metric space $\widehat{A} \triangleq \widehat{C_{\mathbf{s}}[0_{\mathbf{s}}, 1_{\mathbf{s}}]_{\mathbf{s}}}$, where $C_{\mathbf{s}}[0_{\mathbf{s}}, 1_{\mathbf{s}}]_{\mathbf{s}}$ is a paraconsistent linear space paraconsistent **s**-continuous functions $f_{\mathbf{s}} : [0_{\mathbf{s}}, 1_{\mathbf{s}}]_{\mathbf{s}} \to_{\mathbf{s}} \mathbb{R}^{\#}$.
(**b**) With the $w$-metric

$$d_w(f,g) =_w w\text{-}\sup_{0_w \leq_w x \leq_w 1_w} |f_1(x) - f_2(x)|_w \tag{3.2.41}$$

we have $w$-paracomplete separable paraconsistent metric space $\widehat{A} \triangleq \widehat{C_w[0_w, 1_w]_w}$, where $C_w[0_w, 1_w]_w$ is a paraconsistent linear space paraconsistent $w$-continuous functions $f_w : [0_w, 1_w]_w \to_w \mathbb{R}^{\#}$.

# IV.Berry's and Richard's inconsistent numbers within $\mathbf{Z}_2^{\#}$.

## IV.1.Hierarchy Berry's inconsistent numbers $\mathbf{B}_n^{w,(m)}$.

Suppose that $\mathcal{F}_{n_1}^{w}(n,X) \in L_2^{\#}$ is a well-formed formula of second-order arithmetic $\mathbf{Z}_2^{\#}$, i.e. formula which is arithmetical, which has one free set variable $X$ and one free individual variable $n$. Suppose that $g(\exists X \mathcal{F}_{n_1}^{w}(x,X)) \leq \mathbf{k}$, where $g(\exists X \mathcal{F}_{n_1}^{w}(x,X))$ is a corresponding Gödel number. Let be $A_{\mathbf{k}}^{w}, \mathbf{k} \in \mathbb{N}$ the set of all positive weakly inconsistent integers $\bar{n} \in_{\mathbf{s}} \mathbb{N}_w^{\#}$ which can be defined within $\mathbf{Z}_2^{\#}$ (in weak inconsistent sense) under corresponding well-formed formula $\mathcal{F}_{\bar{n}_1(\bar{n})}(x,X)$, i.e. $\exists X_{\bar{n}} \forall m \left[\mathcal{F}_{\bar{n}_1(\bar{m})}^{w}(\bar{m}, X_{\bar{n}}) \to \bar{m} =_w \bar{n}\right]$, hence $\bar{n} \in A_{\mathbf{k}}^{w} \leftrightarrow \exists X_{\bar{n}} \mathcal{F}_{\bar{n}_1(\bar{n})}^{w}(\bar{n}, X_{\bar{n}})$.

Thus

$$\forall n_{n\in_{\mathbf{s}}\mathbb{N}_w^\#}[n \in A_{\mathbf{k}}^w \leftrightarrow \exists X_n \mathcal{F}_{n_1(n)}(n, X_n)],$$
$$\text{where} \tag{4.1.1}$$
$$g(\exists X_n \mathcal{F}_{n_1}(x, X_n)) \leq \mathbf{k},$$

Since there are only finitely many of these $\bar{n}$, there must be a smallest (relative to $<_w$ ) positive integer $\mathbf{B}_{\mathbf{k}}^w \in_w \mathbb{N}_w^\# \backslash_w A_{\mathbf{k}}^w$ that does not belong to $A_{\mathbf{k}}^w$.But we just defined $\mathbf{B}_{\mathbf{k}}^w$ in under corresponding well-formed formula

$$\mathbf{B}_{\mathbf{k}}^w \in_w A_{\mathbf{k}}^w \leftrightarrow \breve{\mathcal{F}}_{\bar{n}_1(\mathbf{B}_{\mathbf{k}}^w)}(\mathbf{B}_{\mathbf{k}}^w, A_{\mathbf{k}}^w),$$
$$\tag{4.1.2}$$
$$\breve{\mathcal{F}}_{\bar{n}_1(\mathbf{B}_{\mathbf{k}})}(\mathbf{B}_{\mathbf{k}}^w, A_{\mathbf{k}}^w) \leftrightarrow \mathbf{B}_{\mathbf{k}}^w =_w w\text{-}\min_{n\in_w \mathbb{N}_w^\#}(\mathbb{N}_w^\# \backslash_w A_{\mathbf{k}}^w).$$

Hence for a sufficiently Large $\mathbf{k}$ such that: $g\left(\breve{\mathcal{F}}(\mathbf{B}_{\mathbf{k}}^w, A_{\mathbf{k}})\right) \leq \mathbf{k}$ we obtain:

$$(\mathbf{B}_{\mathbf{k}}^w \in_w A_{\mathbf{k}}^w) \wedge (\mathbf{B}_{\mathbf{k}}^w \notin_w A_{\mathbf{k}}^w). \tag{4.1.3}$$

**Theorem.4.1.1**. Paraconsistent set $A_{\mathbf{k}}^w \subset_w \mathbb{N}_w^\#$ wich was defined above it a strictly $\in$-inconsistent set with $\mathbf{rank} \geq 0$.

# IV.2.Hierarchy Richard's inconsistent numbers $\Re_n^{w,(m)}$.

Let be $q_n \in_{\mathbf{s}} \mathbb{Q}_w^\#$ paraconsistent rational number with corresponding decimal representation $q_n =_{\mathbf{s}} \{0, q_n(1_w)q_n(2_w)\ldots q_n(i)\ldots q_n(n)\}, n \in_{\mathbf{s}} \mathbb{N}_w^\#$,

$q_n(i) =_w 0_w \vee 1_w \vee 2_w \vee 3_w \vee 4_w \vee 5_w \vee 6_w \vee 7_w \vee 8_w \vee 9_w$,

$i \leq_w n$, $x_k^w =_{\mathbf{s}} \langle q_n^k \rangle =_{\mathbf{s}} \langle q_n^k : n \in \mathbb{N}, q_n^k =_{\mathbf{s}} \{0, q_n^k(1)q_n^k(2)\ldots q_n^k(n)\}\rangle_{\mathbf{s}} \in_{\mathbf{s}} \mathbb{R}_w^\#, k \in \mathbb{N}$
is a paraconsistent real number which can be defined in weak inconsistent sense under corresponding well-formed formula (of second-order arithmetic $\mathbf{Z}_2^\#$) $\mathcal{F}_k(x, X)$, i.e. $\forall q(q \in_{\mathbf{s}} \mathbb{Q}_w^\#)[q \in_w \langle q_n^k \rangle \leftrightarrow \exists X \mathcal{F}_k(q, X)]$.

**Definition**. **4.2.1**. We denote paraconsistent real number $x_k^w$ as $k$-th Richard's weakly inconsistent real number.

Let us consider Richard's real number $\Re_p^w =_{\mathbf{s}} \langle \Re_n^p : n \in \mathbb{N} \rangle$ such that

$$\mathfrak{R}_n^p =_w 1_w \leftrightarrow q_n^n(n) \neq_w 1_w,$$
$$\mathfrak{R}_n^p =_w 0_w \leftrightarrow q_n^n(n) =_w 1_w. \tag{4.2.1}$$

Suppose that $q_p^p(p) \neq_w 1$, hence $\mathfrak{R}_p^p(p) =_w 1_w$. Thus $\mathfrak{R}_p^p(p) \neq_w q_p^p(p) \to \mathfrak{R}_p \neq_w x_p$.
Suppose that $q_p^p(p) =_w 1_w$, hence $\mathfrak{R}_p^p(p) =_w 0_w$. Thus $\mathfrak{R}_p^p(p) \neq_w q_p^p(p) \to \mathfrak{R}_p \neq_w x_p$.
Hence for any Richard's real number $x_k$ one obtain the contradiction

$$\forall k\, [x_k \neq_w \mathfrak{R}_p^p(p)]. \tag{4.2.2}$$

**Theorem.4.2.1**. Paraconsistent set $\{x_k^w : k \in \mathbb{N}\}_{\mathbf{s}} \subset_w \mathbb{R}_w^\#$ containing the all Richard's weakly inconsistent real numbers which was defined above it a strictly $\in$-inconsistent set with $\mathbf{rank} \geq 0$.

# V.Schinzel's conjecture H within $Z_2^\#$.

## V.1.Schinzel's conjecture H.

Schinzel's conjecture H is a very broad generalisation of conjectures such as the twin prime conjecture. It aims to define the possible scope of a conjecture of the nature that several sequences of the type

$$f_1(n), f_2(n), \ldots, f_k(n) \tag{5.1.1}$$

with values at integers $n$ of irreducible integer-valued polynomials $f_1(x), f_2(x), \ldots, f_k(x)$ should be able to take on prime number values simultaneously, for integers $n$ that can be as large as we please. Putting it another way, there should be infinitely many such $n$, for which each of the sequence values are prime numbers. Some constraints are needed on the polynomials. Andrzej Schinzel's hypothesis builds on the earlier Bunyakovsky conjecture, for a single polynomial, and on the Hardy-Littlewood conjectures for multiple linear polynomials. It is in turned extended by the Bateman-Horn conjecture.

The arithmetic nature of the most evident necessary conditions can be understood. An integer-valued polynomial $Q(x)$ has a fixed divisor $m$ if there is an integer $m > 1$ such that

$$\frac{Q(x)}{m} \tag{5.1.2}$$

is also an integer-valued polynomial. Such fixed divisors must be ruled out of

$$Q_k(x) = \prod_{i=1}^{k} f_i(x) \tag{5.1.3}$$

for any conjecture for polynomials $f_i, i = 1, \ldots, k,$ since their presence is easy seen to contradict the possibility that $f_i(n)$ can all be prime, with large values of $n$.

Therefore the standard form of hypothesis H is that if $Q_k(x)$ defined as above has no fixed prime divisor, then all $f_i(n)$ will be simultaneously prime, infinitely often, for any choice of irreducible integral polynomials $f_i(x)$ with positive leading coefficients.If the leading coefficients were negative, we could expect negative prime values; this is a harmless restriction, really. There is probably no real reason to restrict to integral polynomials, rather than integer-valued polynomials. The condition of having no fixed prime divisor is certainly effectively checkable in a given case, since there is an explicit basis for the integer-valued polynomials. As a famous simple example:

$$x^2 + 1 \tag{5.1.4}$$

has no fixed prime divisor. We therefore expect that there are infinitely many primes:

$$n^2 + 1. \tag{5.1.5}$$

This has not been proved, though. It was one of Landau's conjectures and goes back to Euler, who observed in a letter to Goldbach in 1752 that $n^2 + 1$ is often prime for $n$

up to 1500.

At the 1912 International Congress of Mathematicians, Edmund Landau listed four basic problems about primes. These problems were characterised in his speech as "unattackable at the present state of science" and are now known as Landau's problems. They are as follows:

**1**.Goldbach's conjecture: Can every even integer greater than 2 be written as the sum of two primes?

**2**.Twin prime conjecture: Are there infinitely many primes $p$ such that $p+2$ is prime?

**3**.Legendre's conjecture: Does there always exist at least one prime between consecutive perfect squares?

**4**.Are there infinitely many primes $p$ such that $p-1$ is a perfect square? In other words: Are there infinitely many primes of the form $n^2+1$?

**Theorem 5.1.1**

## V.2. Almost prime natural numbers.

Let $\Omega(n)$ be the total number of primes in the prime factorization of $n,$ i.e.

$$\Omega(n) = \sum_{i=1}^{r} m_i \leftrightarrow n = \prod_{i=1}^{r} p_i^{m_i}. \qquad (5.2.1)$$

**Definition 5.2.1**. [8][see [8] p.48] An natural number $n$ is called $k$-almost prime if it has exactly $k$ prime factors, counted with multiplicity,i.e. a number $n$ is $k$-almost prime if and only if $\Omega(n) = k$.

$$(5.2.)$$

## V.3. Some General Results on Dirichlet Series over paraconsistent field $\mathbb{C}^{\#}_{\mathbf{inc}}$.

**Definition 5.3.1**. A Dirichlet $w$-series over paraconsistent field $\mathbb{C}^{\#}_{\mathbf{inc}}$ is a $w$-series of the form:

$$F_w(a_n,s) = \sum_{n=1}^{\infty} a_n \times_w n^{-s}, s \in \mathbb{C}, \qquad (5.3.1)$$

where $\{a_n\}_{n\in\mathbb{N}}$ is some $w$-sequence of $w$-complex numbers such that $a_n \in_{\mathbf{s}} \mathbb{C}^{\#}_{\mathbf{inc}}$.

We will want to know when a Dirichlet $w$-series over paraconsistent field $\mathbb{C}^{\#}_{\mathbf{inc}}$ actually $w$-converges. In order to prove this, we will need the following simple theorem which is the analog for sums of $w$-integration by parts.

**Theorem 5.3.1**. (Generalized Abel Summation Formula or Summation by Parts). Suppose $(a_n)$ and $(b_n)$ are two $w$-sequences and $N \in \mathbb{N}$. Then

$$w\text{-}\sum_{n=l}^{N} a_n \times_w (b_n -_w b_{n-1}) = a_N \times_w b_N -_w a_l \times_w b_l - w\text{-}\sum_{n=l}^{N} b_n \times_w (a_n - a_{n-1}) \qquad (5.3.2)$$

Proof. To verify this theorem, one simply check that each term $a_i \times_w b_j$ occurs with the same multiplicity on both sides.Using the notion of a Stieltjes $w$-integral we can rewrite the Abel Summation Formula as a generalization of $w$-integration by parts.

**Definition 5.3.2**. Given a $w$-closed $w$-interval $I_w = [a,b]_w$, a $w$-partition of $I_w$ is any finite strictly $w$-increasing $w$-sequence of points $P_w =_w \{x_0,x_1,\ldots,x_{n-1},x_n\}_w$such that $a =_w x_0$ and $b =_w x_n$. Each partition $P_w =_w \{x_0,x_1,\ldots,x_{n-1},x_n\}_w$ of $I_w$ decomposes $I_w$ into $n$ $w$-subintervals $I_{j,w} =_w [x_{j-1,w},x_{j,w}]_w, j = 1,\ldots,n,$ such that $I_{j,w} \cap_w I_{k,w} =_w \varnothing$ for $k \neq j$ or $k \neq j+1$. The set of the all partitions as above we denoted by $\Pi_w(a,b)$.

**Definition 5.3.3**. Let $\alpha$ be a non-$w$-decreasing (not necessarily $w$-continuous $w$-function) on $[a,b]_w$ and $f$ is a $w$-function $w$-bounded on $[a,b]_w$. Let $P_w$ be some $w$-partition of $[a,b]_w$,let $\Delta_w\alpha_i =_w \alpha_i -_w \alpha_{i-1}$. Define the $w$-upper and $w$-lower $w$-sums, $U_w = U_w(P_w,f,\alpha) = w\text{-}\sum_{i=1}^{n} M_{w,i} \times_w \alpha_i$ and $L_w = L_w(P_w,f,\alpha) = w\text{-}\sum_{i=1}^{n} m_{w,i} \times_w \alpha_i$where $M_{w,i}$ and $m_{w,i}$ are the $w$-maximum and $w$-minimum respectively of $f$ on the $i$-th $w$-interval in the $w$-partition. We define the $w$-upper and $w$-lower Stieltjes integrals to be

$$w\text{-}\overline{\int_a^b} f(x)d\alpha(x) = w\text{-}\inf_{P_w \in \Pi_w(a,b)} U_w(P_w, f, \alpha) \tag{5.3.3}$$

and

$$w\text{-}\underline{\int_a^b} f(x)d\alpha(x) = w\text{-}\sup_{P_w \in \Pi_w(a,b)} L_w(P_w, f, \alpha) \tag{5.3.4}$$

respectively. If

$$w\text{-}\overline{\int_a^b} f(x)d\alpha(x) =_w w\text{-}\underline{\int_a^b} f(x)d\alpha(x). \tag{5.3.5}$$

Then $f$ is Riemann-Stieltjes $w$-integrable, or $w$-integrable with respect to $\alpha(x)$ in the Riemann sense, on $[a,b]_w$ , and the common value of the integral is denoted by

$$I_w[f;a,b] \triangleq w\text{-}\int_a^b f(x)d\alpha(x) \tag{5.3.6}$$

and say that $I_w[f;a,d]$ is $w$-integral with respect to $\alpha$.For our purposes we will only be looking at Stieltjes $w$-integrals where either $\alpha$ is $w$-differentiable in which case

$$w\text{-}\int_a^b f(x)d\alpha(x) =_w w\text{-}\int_a^b f(x)[w\text{-}\alpha'(x)]dx. \tag{5.3.7}$$

or when $\alpha$ is a step function. In the latter situation the $w$-integral is simply a $w$-sum over the jumps of the value of the function at that jump times the size of the jump.

**Proposition 5.3.1.**If $w\text{-}\sum_{n=1}^{\infty} a_n \times_w n^{-s_0}, s \in \mathbb{C}^{\#}_{\mathbf{inc}}$ $w$-converges, then the Dirichlet $w$-series $F_w(a_n,s)$ $w$-converges for all complex numbers $s \in \mathbb{C}^{\#}_{\mathbf{inc}}$ with $\mathrm{Re}(s)\ _w > s_0$. In fact, this $w$-convergence is uniform in any wedge to the right of the point $s_0$ : $\left\{s : \mathrm{Re}(s)_w > \mathrm{Re}(s_0) \text{ and } 0 <_w \frac{|s -_w s_0|_w}{\mathrm{Re}(s -_w s_0)} <_w M\right\}_w$ , where $M$ is an arbitrary $w$-positive constant. (Since this $w$-convergence is uniform, $F_s(a_n)$ is $w$-analytic on that

region.)

**Proposition 5.3.2**. Suppose that $w\text{-}\sum_{n=1}^{\infty} a_n \times_w n^{-s} =_w 0_w$ on some right half plane $\mathrm{Re}(s) > \sigma_0$.Then, $a_n = 0_w$ for all $n \in \mathbb{N}$. Therefore, if we have two Dirichlet $w$-series with $F_w(a_n,s) =_w F_w(b_n,s)$ for all $s : \mathrm{Re}(s)_w > \sigma_0$, then $a_n =_w b_n$ for all $n \in \mathbb{N}$.

**Definition 5.3.4**. An $w$-sequence $a_n$ is called $w$-multiplicative if $a_n \times_w a_m =_w a_{nm}$ for all relatively prime positive integers $n \in \mathbb{N}$ and $m \in \mathbb{N}$. Similarly, an $w$-sequence $a_n$ is called strongly $w$-multiplicative if $a_n \times_w a_m =_w a_{nm}$ for all pairs of positive integers.

**Theorem 5.3.2**. If $a_n$ is $w$-multiplicative then the Dirichlet $w$-series $F_w(a_n,s)$ has the Euler factorization such that:

$$F_w(a_n,s) =_w w\text{-}\prod_p \left( w\text{-}\sum_{k=1}^{\infty} a_{p^k} \times_w p^{-sk} \right). \tag{5.3.8}$$

Furthermore if $a_n$ is strongly $w$-multiplicative, summing this geometric $w$-series we have that

$$F_w(a_n,s) =_w w\text{-}\prod_p \left(1 - \frac{a_p}{p^s}\right)^{-1_w}. \tag{5.3.9}$$

Finally we have a formula for the $w$-product of two Dirichlet $w$-series.

$$F_w(a_n,s) \times_w F_w(b_n,s) =_w F_w\left(w\text{-}\sum_{d|n} a_d b_{\frac{n}{d}}, s\right). \tag{5.3.10}$$

**Proof**. By definition one obtain

$$F_w(a_n,s) F_w(b_n,s) =_w w\text{-}\sum_{m=1}^{\infty} w\text{-}\sum_{l=1}^{\infty} \frac{a_m \times_w b_l}{(ml)^n} \tag{5.3.11}$$

Make a change of variables $n = ml$ and $d = m$, one obtain

$$F_w(a_n,s) \times_w F_w(b_n,s) =_w w\text{-}\sum_{n=1}^{\infty}\left(w\text{-}\sum_{d|n}\frac{a_d b_{\frac{n}{d}}}{n^s}\right) =_w F_w\left(w\text{-}\sum_{d|n} a_d b_{\frac{n}{d}}, s\right). \qquad (5.3.12)$$

# V.4. Dirichlet's $L^{\#}$-series over paraconsistent field $\mathbb{C}^{\#}_{\mathbf{inc}}$.

In order to generalize Dirichlet's argument from the case $lx+m$ to a general polynomials $f(x)$ we should look at classical Dirichlet series. From our results above, it is clear that we should be looking at the $w$-series corresponding to $w$-sequences of the following form:

**Definition 5.4.1**. Let $G$ be a finite classical, i.e. consistent abelian group, written multiplicatively and with unity $e$. A $w$-character on $G$ is a $w$-homomorphism

$$\chi_w : G \to (\mathbb{C}^{\#}_{\mathbf{inc}})^{\times}. \qquad (5.4.1)$$

The image of $\chi_w$ is contained in the $w$-group of $|G|$-th $w$-roots of $w$-unity $1_w$.The $w$-characters form a group themselves, called the dual $w$-group of $G$ and written $G^{*}_{\mathbf{inc}}$. The identity element in this group is $\chi_{w;0}$, the trivial or principal $w$-character which maps everything in $G$ to $w$-unity $1_w$.

**Definition 5.4.2**. Let a Dirichlet $w$-character modulo $m \in \mathbb{N}$ be any $w$-function $\chi_{m,w} : \mathbb{Z} \mapsto \mathbb{C}^{\#}_{\mathbf{inc}}$ with the properties:

**1**. If $n$ and $m$ are not relatively prime, then $\chi_{m,w}(n) =_w 0_w$;

**2**. If $n$ and $m$ are relatively prime, then $|\chi_{m,w}(n)|_w =_w 1_w$;

**3**. If $n_1$ and $n_2$ are any two positive integers, then $\chi_{m,w}(n_1 \times n_2) =_w \chi_{m,w}(n_1) \times_w \chi_{m,w}(n_2)$.

**Definition 5.4.3**.If $\chi_w$ is a Dirichlet $w$-character modulo $m,$ then define the Dirichlet $L^{\#}$-series

$$L^{\#}(s,\chi_w) =_w w\text{-}\sum_{n=1}^{\infty}\chi_w(n) \times_w n^{-_w s} =_w w\text{-}\prod_p \frac{1_w}{\left(1_w -_w \frac{\chi_w(p)}{p^s}\right)}. \qquad (5.4.2)$$

If $\chi^0_w$ is the trivial character modulo $m,$ then obviously

$$L^{\#}(s,\chi_w^0) =_w w\text{-}\prod_{p|m}\left(1_w -_w \frac{\chi_w^0(p)}{p^s}\right) \times_w \zeta_w(s),$$

$$\zeta_w(s) =_w w\text{-}\sum_{n=1}^{\infty} n^{-_w s} \tag{5.4.3}$$

Dirichlet $L^{\#}$-series $L^{\#}(s,\chi_w^0)$ $w$-converges for any $\mathrm{Re}(s)\ _w > 1_w$.

**Lemma 5.4.1.**(Generalized Dirichlet Lemma) Let $G$ be a finite classical, i.e. consistent abelian group, written multiplicatively and $\chi_w$ is a nontrivial $w$-character on $G$, then

$$w\text{-}\sum_{g\in G} \chi_w(g) =_w 0_w. \tag{5.4.4}$$

**Theorem 5.4.1.**(Generalized Dirichlet Theorem) (1) Let $G$ be a finite classical, i.e. consistent abelian group, written multiplicatively. Let $\chi_{w,1}$ and $\chi_{w,2}$ be a nontrivial $w$-character on $G$, then

$$w\text{-}\sum_{g\in G} \chi_{w,1}(g) \times_w \chi_{w,2}^{-1}(g) =_w 0_w. \tag{5.4.5}$$

(2) Let $G_w^*$ be the $w$-dual group of $w$-characters $\chi_w$ of $G$, then, if $g \neq_w h$

$$\sum_{\chi_w\in G_w^*} \chi_w(gh^{-1}) =_w 0_w. \tag{5.4.6}$$

**Remark 5.4.1.**Thus by Theorem 5.4.1, we see that the $w$-characters $\chi_w$ are all $w$-linearly independent and thus form an $w$-basis of the space of all complex $\mathbb{C}_{\mathbf{inc}}^{\#}$-valued functions on $G$. In particular,

$$\frac{1}{|G|}\times_w\left(w\text{-}\sum_{\chi_w\in G_w^*}\chi_w(gh^{-1})\right)=_w\begin{cases}1_w \text{ if } g=h\\ 0_w \text{ otherwise}\end{cases} \tag{5.4.7}$$

Therefore, if we let $G=(\mathbb{Z}/m\mathbb{Z})^\times$ then for $(l,m)=1$ one obtain

$$\frac{1}{\varphi(m)}\times_w\left(w\text{-}\sum_{\chi_w}\chi_w(n)\times_w\overline{\chi}_w(l)\right)=_w\begin{cases}1_w \text{ if } n\equiv l \bmod m\\ 0_w \text{ if } n\not\equiv l \bmod m\end{cases} \tag{5.4.8}$$

We notice that

$$w\text{-}\log L^\#(s,\chi_w)=_w w\text{-}\sum_p w\text{-}\log\frac{1_w}{\left(1_w-_w\frac{\chi_w(p)}{p^s}\right)}. \tag{5.4.9}$$

Obviously the Taylor $w$-series expansion is valid, and our estimate on the terms still holds:

$$w\text{-}\sum_{n=2}^{\infty}\frac{1_w}{n}\times_w\left(w\text{-}\sum_p\frac{\chi_w(p)}{p^{n\times_w s}}\right)<_w w\text{-}\int_2^\infty w\text{-}\int_2^\infty x^{-_w 1_w}y^{-_w n\times_w s}dydx<_w 32^{-_w 1_w}. \tag{5.4.10}$$

Therefore,

$$w\text{-}\log L^\#(s,\chi_w)=_w w\text{-}\sum_p\chi_w(p)\times_w p^{-_w s}+_w O(1_w). \tag{5.4.11}$$

**Definition 5.4.4**. Let a Dirichlet $w_{[1]}$-character modulo $m\in\mathbb{N}$ coresponding to polynomial $x^2+1$, be a $w_{[1]}$-function $\chi_{m,w,[1]}:\mathbb{Z}\mapsto\mathbb{C}^\#_{\mathbf{inc}}$ with the properties:

**1**. If $n$ and $m$ are not relatively prime, then $\chi_{m,w,[1]}(n) =_{w,[1]} 0_w$;
**2**. If $n$ and $m$ are relatively prime, then $|\chi_{m,w,[1]}(n)|_w =_{w,[1]} 1_w$;
**3**. If $n_1$ and $n_2$ are any two positive integers, then
$\chi_{m,w,[1]}(n_1 \times n_2) =_{w,[1]} \chi_{m,w,[1]}(n_1) \times_w \chi_{m,w,[1]}(n_2)$;
**4**. For $(l,m) = 1$

$$\frac{1}{\varphi(m)} \times_w \left( w\text{-}\sum_{\chi_w} \chi_{w,[1]}(n) \times_w \overline{\chi}_{w,[1]}(l) \right) =_{w,[1]}$$

$$=_{w,[1]} \begin{cases} 1_w \text{ if } (n \equiv l \bmod m) \wedge (n = r^2+1) \\ 0_w \text{ if } (n \not\equiv l \bmod m) \vee (n \neq r^2+1) \end{cases} \quad (5.4.12)$$

In order to get main result,we need to write

$$f_{l,m,w}(s) =_w w\text{-}\sum_{(p \equiv l \bmod m) \wedge (p = r^2+1)} p^{-_w s} \quad (5.4.14)$$

as a sum of $w\text{-}\log L^{\#}(s, \chi_{w,[1]})$ for various $\chi_{w,[1]}$.

Since $w\text{-}\log L^{\#}(s, \chi_{0,w,[1]})$ blows up, and $w\text{-}\log L^{\#}(1/2, \chi_{0,w,[1]}) <_{w,[1]} \infty$ for all nontrivial $w_{[1]}$-characters. All that remains is to show that $w\text{-}\log L^{\#}(1/2, \chi_{0,w,[1]})\ _{w,[1]} > -_{w,[1]} \infty$.

**Theorem 5.4.2.**If $\chi_{w,[1]}$ is a nontrivial Dirichlet $w_{[1]}$-character coresponding to polynomial $x^2+1$, then $L^{\#}(1/2, \chi_{0,w,[1]}) \neq_{w,[1]} 0_{w,[1]}$.